\documentclass[intlimits,12pt,psamsfonts]{amsart}
\newtheorem{thm}{Theorem}
\newtheorem*{thm*}{Theorem}
\newtheorem{lemma}[thm]{Lemma}
\newtheorem{cor}[thm]{Corollary}
\newtheorem{prop}[thm]{Proposition}
\theoremstyle{definition}
\newtheorem*{defn}{Definition}
\newtheorem*{example}{Example}
\theoremstyle{remark}
\newtheorem*{rmk}{Remark}
\newcommand{\leg}[2]{\genfrac{(}{)}{}{}{#1}{#2}}
\newcommand{\re}{\mathrm{Re}}
\newcommand{\im}{\mathrm{Im}}
\newcommand{\Li}{\mathrm{Li}}
\usepackage{graphicx}
\usepackage{amssymb}
\usepackage{psfrag}
\usepackage{flafter}
\begin{document}
\title{Hexagonal Lattice Points on Circles}
\author{Oscar Marmon}
\begin{abstract} 

We study the angular distribution of points in the hexagonal lattice
(i.e., $\mathbb{Z}[\frac{1+\sqrt{-3}}{2}]$) lying on a circle centered
at the origin.  We prove that the angles are equidistributed on
average, and show that the discrepancy is quite small for almost all
circles.  Equidistribution on average is expressed in terms of
cancellation in exponential sums. We introduce Hecke L-functions and
investigate their analytic properties in order to derive estimates on
sums of Hecke characters. Using a version of the Halberstam-Richert
inequality, these estimates then yield the desired results for the
exponential sums.

An interesting consequence of these bounds is that the discrete velocity model
(DVM) for the Boltzmann equation is consistent when using a hexagonal
lattice.

\end{abstract}
%-----------------------------------------------------------------------
\maketitle
%-----------------------------------------------------------------------
\newpage
\tableofcontents
%-----------------------------------------------------------------------
\newpage
\section{Introduction}\label{intro}
%-----------------------------------------------------------------------
\subsection{The Hexagonal Lattice}
\begin{figure}
\centering
\includegraphics[%trim=30 0 0 0,
scale=0.75,clip]{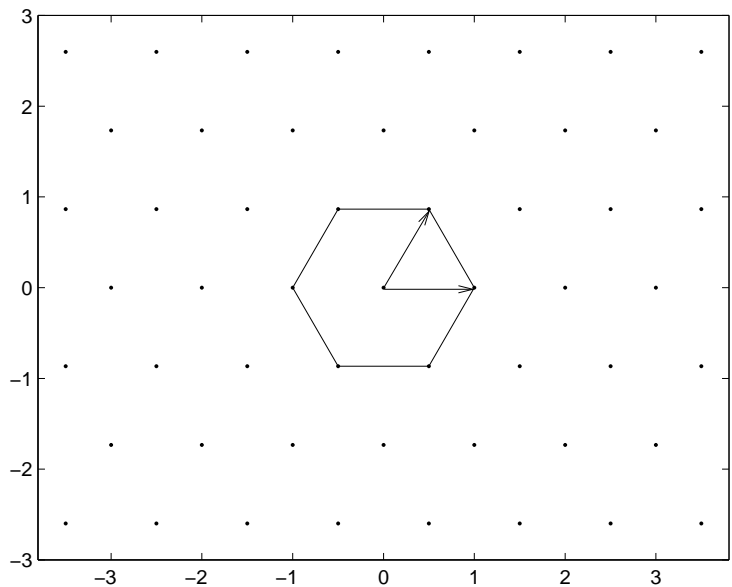}
\caption{The hexagonal lattice}
\label{lattice}
\end{figure}

We will study the properties of a lattice in $\mathbb{R}^2$ that is spanned by the vectors $(1,0)$ and $\bigl(\frac{1}{2},\frac{\sqrt{3}}{2}\bigr)$,
see Figure \ref{lattice}.
We are interested in the following two questions: given a circle with radius $r$, centered at the origin, whose perimeter contains at least one lattice point, \begin{enumerate}\label{questions}
\item\label{q1} how many lattice points lie on the circle and
\item how are these distributed around the circle?
\end{enumerate}
\begin{figure}
\centering
\includegraphics[scale=0.75]{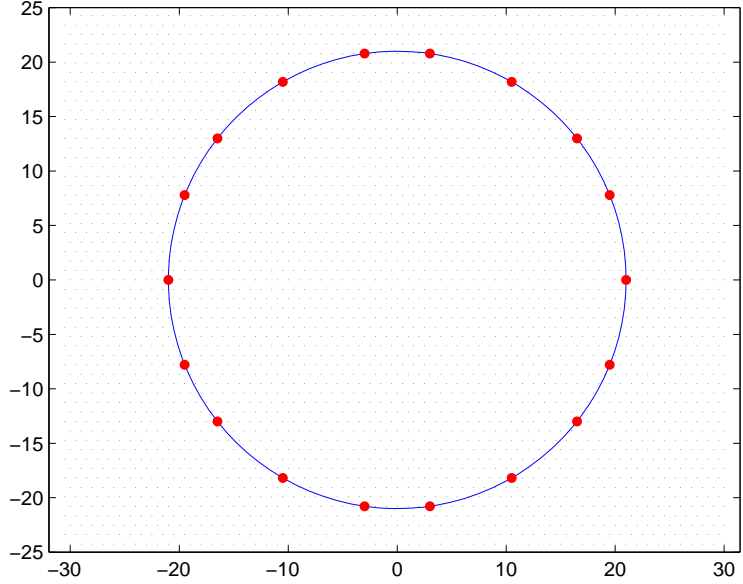}
\caption{Lattice points on the circle with radius $r=21$}
\label{circle21}
\end{figure}
Why are these questions interesting? Well, aside from the intrinsic number theoretic and geometric importance, there is, perhaps surprisingly, an application to the Boltzmann equation in  the kinetic theory of gases.
%-----------------------------------------------------------------------
\subsection{The Boltzmann Equation}
Under certain hypotheses, the behavior of a gas is described by the \emph{Boltzmann Equation}:
\begin{equation}\label{boltzmann}
\frac{\partial{f}}{\partial{t}} + v\frac{\partial{f}}{\partial{x}} = Q(f,f).
\end{equation}
Here $f(x,v, t):\mathbb{R}^3\times\mathbb{R}^3\times\mathbb{R}^+\to\mathbb{R}^+$ is the \emph{phase space density} of the gas, that is, it describes the expected mass density of the gas at time $t$ in the point $(x,v)$ in six-dimensional phase-space, where $x\in\mathbb{R}^3$ is a position in space and $v\in\mathbb{R}^3$ is a velocity. \\
If we replace $3$ with any spatial dimension $d$, the Boltzmann equation, although losing its physical sense, is still of mathematical interest. We will consider the case $d=2$.  $Q(f,f)$ is called the \emph{collision term} and is given in the 2-dimensional case by
\begin{equation}\label{collision}
Q(f,f)(v)=\int_{\mathbb{R}^2}\left(\int_{-\pi}^\pi\big(f(v')f(v'_*)-f(v)f(v_*)\big)q(|w|,\cos\theta)d\theta\right) dv_*,
\end{equation}
where $v,v_*$ are the velocities after the collision and $v',v'_*$ those before. Moreover, $q(|w|,\cos\theta)$ is a quantity describing the probability that two particles with relative velocity $2w$, i.e.
 \[
 w=\frac{v-v_*}2,
 \]
 collide with relative deflection angle $\theta$. We have the relations
 \begin{equation}\begin{split}
 \begin{cases}v'=\frac12(v+v_*)+|w|u,\\
 v'_*=\frac12(v+v_*)-|w|u,
 \end{cases}
 \end{split}\end{equation}
where $u$ is a unit vector with $u=(\cos\theta,\sin\theta)$, which follow from the laws of conservation of energy and momentum. In other words, $v'$ and $v'_*$ are endpoints of a diameter of a circle with radius $|w|$, centered at $(v+v_*)/2$.\\  
In a \emph{discrete velocity model} (DVM), one considers a discrete set of possible velocities. (See \cite{bobylev} for more details.) Here we will consider a DVM where the velocities belong to the set $h\mathbb{L}$, where $h>0$ and $\mathbb{L}$ is some lattice in $\mathbb{R}^2$. Set
\[
f^h=\sum_{\xi\in\mathbb{L}}f\delta_{v=h\xi},
\]
so that $f^h\to f$ as $h\to 0$, in some suitable sense. We want to prove the \emph{consistency} of the DVM, i.e. the property that
\begin{equation}\label{consistency}
Q(f^h,f^h)(v)\to Q(f,f)(v) \text{ for all }v\in h\mathbb{L},\text{ as } h\to 0.
\end{equation} 
Let
\[
g_v(w,\theta)=\big(f(v')f(v'_*)-f(v)f(v_*)\big)q(|w|,\cos\theta).
\]
Then, by a change of variables,
\begin{equation}\label{qwithg}
Q(f,f)(v)=4\int_{\mathbb{R}^2}\left(\int_{-\pi}^\pi g_v(w,\theta)d\theta\right)dw.
\end{equation}
When we discretize the integral in \eqref{qwithg}, the outer integral turns into a sum over all lattice points $w=h\zeta$, $\zeta\in\mathbb{L}$, while for the inner integral we get a sum over all $u$ such that $v'$ and $v'_*$ belong to $h\mathbb{L}$. A necessary and sufficient condition for this is that $u=\frac{\xi}{|\xi|}$, where $\xi\in \mathbb{L}$ and $|\xi|=|\zeta|$. Thus we get 
\[
Q(f^h,f^h)(v)=h^2T_{\mathbb{L}}\sum_{\zeta\in\mathbb{L}}\left(\frac{1}{r_{\mathbb{L}}(|\zeta|^2)}\sum_{\substack{\xi\in\mathbb{L}\\|\xi|=|\zeta|}}g_v(h\zeta,\arg\xi)\right),
\]
where $T_{\mathbb{L}}$ is the area of a fundamental region of $\mathbb{L}$ and\[
r_\mathbb{L}(n):=\#\{\xi\in\mathbb{L};|\xi|^2=n\}.
\]

It is clear that $g_v(w,\theta)$ is a $2\pi$-periodic function of $\theta$, and so can be expressed as a Fourier series. It turns out that, assuming certain regularity conditions on $g_v(w,\theta)$, a sufficient condition for \eqref{consistency} to hold is that we have sufficient nontrivial cancellation in the exponential sums  
\[
S(n,A):=\sum_{\substack{\mu\in\mathbb{L}\\|\mu|^2=n}}e^{iA\arg(\mu)}
\]
as $n$ grows, if $A\in\mathbb{Z}$, $A\neq 0$.\\
Fainsilber, Kurlberg and Wennberg proved \cite{kurlberg} that if we choose $\mathbb{L}$ to be $\mathbb{Z}^2$, then \eqref{consistency} holds provided $g_v(w,\theta)$ is a $C^2$-function. This is done by identifying $\mathbb{Z}^2$ with the ring of Gaussian integers, $\mathbb{Z}[i]\subset\mathbb{C}$, and applying number theoretic methods. The results of this paper indicate that we can equally well consider the hexagonal lattice defined above, corresponding to the ring of integers in the algebraic number field $\mathbb{Q}(\sqrt{-3})$. Similar methods to the ones used here should give the same results for every imaginary quadratic number field $\mathbb{Q}(\sqrt{-d})$ with class number 1 (and thus unique factorization into irreducibles). There are exactly 9 such fields: $\mathbb{Q}(\sqrt{-d})$ has unique factorization if and only if $d$ takes one of the values 1, 2, 3, 7, 11, 19, 43, 67 or 163 (see \cite{ireland}, Ch.13,\S 1).\\
Moreover, it is reasonable to expect \eqref{consistency} to hold for any lattice $\mathbb{L}$ such that
\[
\limsup_{n\to\infty}r_{\mathbb{L}}(n)=\infty.
\]
For dimensions $d\geq 3$ (i.e., including the physically relevant case $d=3$), the consistency of a DVM based on the lattice $\mathbb{Z}^d$ was proven by Bobylev, Palczewski and Schneider \cite{bobylev}, using deep number theoretic results. 
%-----------------------------------------------------------------------
\subsection{A Number Theoretic Point of View}\label{subsectNT}
To answer the questions on Page \pageref{q1}, we will view the lattice not as a subset of $\mathbb{R}^2$ but of $\mathbb{C}$, and use number theoretic methods. We identify the lattice with the set
\[
\mathbb{Z}[\omega]=\{a+b\omega\ ;\ a,b\in\mathbb{Z}\},
\]
where
\[
\omega=e^{\frac{i\pi}{3}}=\frac{1}{2}+\frac{i\sqrt{3}}{2}.
\]
This is in fact the ring of integers $\mathfrak{O}$ in the algebraic number field $K:=\mathbb{Q}(\sqrt{-3})$. Since $K$ is a quadratic number field, the norm $N(\alpha)$ of an element $\alpha\in K$ is just the squared complex modulus $|\alpha|^2$. Hence question \eqref{q1} above can be rephrased as:
\begin{quote}
How many elements of $\mathfrak{O}$ are there with a given norm?
\end{quote}
Let $a,b\in\mathbb{Z}$, so that $a+b\omega\in\mathfrak{O}$. Then we have
\[
N(a+b\omega)=(a+b\omega)(a+b\overline{\omega})=a^2+(\omega+\overline{\omega})ab+\omega\overline{\omega}b^2=a^2+ab+b^2,
\]
so that yet another formulation of the first question would be:
\begin{quote}
How many integer solutions $(a,b)$ exist to the equation $a^2+ab+b^2=R$ for a given $R\in\mathbb{Z}$?
\end{quote}
This is a classical question, and the answer we present is by no means new (see for example \cite{hardy} , Ch.~12.9). The norm of a quadratic number field induces a quadratic form, in our case
\[
Q(x,y)=x^2+xy+y^2.
\]
The quantity we seek is then
\[
r_Q(n):=\#\{(a,b)\in\mathbb{Z}^2;Q(a,b)=n\},
\]
the number of representations of $n$ by the form $Q$.  
The key point is, that each lattice point $a+b\omega$ on the circle with radius $r$ corresponds to a representaion of $r^2$ as a product of two elements $a+b\omega$ and $a+b\overline{\omega}$ in $\mathfrak{O}$. Fortunately, $\mathfrak{O}$ is a principal ideal domain, hence a unique factorization domain, so all such representations can easily be found by factorizing $r^2$ into rational primes.\\
For every rational prime $p$, one of the following three cases is true:
\begin{itemize}
\item $p=\pi\cdot\overline{\pi}$ where $\pi$ is a prime in $\mathfrak{O}$. We say that $p$ is \emph{split} in $\mathfrak{O}$.
\item $p=\pi^2$ where $\pi$ is a prime in $\mathfrak{O}$. We say that $p$ is \emph{ramified} in $\mathfrak{O}$.
\item $p$ remains prime in $\mathfrak{O}$. We say that $p$ is \emph{inert} in $\mathfrak{O}$.
\end{itemize}
These factorizations are unique up to multiplication by a unit. The unit group of $\mathbb{Z}[\omega]$ is 
\[
U=\{\pm 1,\pm\omega,\pm\omega^2\}=\{\omega,\omega^2,\ldots,\omega^6\}.
\]
Moreover, all primes in $\mathbb{Z}[\omega]$ can be found in this way. For odd $p$ it is known that (see for example \cite{ireland}, Ch.~13,\S 1) $p$ is split if $\leg{-3}{p}=1$, ramified if $\leg{-3}{p}=0$, and inert if $\leg{-3}{p}=-1$.
Moreover it is easily seen that $2$ is inert, since the equation $a^2+ab+b^2=2$ has no integer solutions. We thus have:
%-----------------------------------------------------------------------
\begin{prop}
Let $p$ be a rational prime. Then
\begin{enumerate}
\item $p$ is split iff $p\equiv 1\pmod{3}$.
\item $p$ is inert iff $p\equiv 2\pmod{3}$.
\item $p$ is ramified iff $p=3$.
\end{enumerate}
\end{prop}
%-----------------------------------------------------------------------
We see that to each rational prime $p$ with $p\equiv 1\pmod{3}$ corresponds a unique prime $\pi_p=\sqrt{p}e^{i\theta_p}$ in $\mathbb{Z}[\omega]$, if we choose $\theta_p$ to lie in the interval $[0,\frac{\pi}{6}]$. Indeed, for every $\pi$, one of the numbers $\omega\pi,\ldots,\omega^6\pi,\omega\overline{\pi},\ldots,\omega^6\overline{\pi}$ must lie in this interval. For p=3, analogously, if we set 
\[
\pi_3=\frac32+i\frac{\sqrt{3}}2=\sqrt{3}e^{i\frac\pi 6},
\]
we have $3=\pi_3\overline{\pi}_3$. We note however that $\overline{\pi}_3=\overline{\omega}\pi_3$, so that $3=\overline{\omega}\pi_3^2$, confirming that $3$ is ramified.\\
Now let $n\in\mathbb{Z}$ and suppose that $n$ has the rational prime factorization
\[
n=3^\alpha p_1^{\alpha_1}\cdots p_k^{\alpha_k}q_1^{\beta_1}\cdots q_l^{\beta_l},\]
where $p_1,\ldots,p_k\equiv 1\pmod{3}$ and $q_1,\ldots,q_l\equiv 2\pmod{3}$.
For each $p_i$ we choose the unique prime $\pi_i=\sqrt{p}e^{i\theta_{p_i}}$. Then the complete factorization of $n$ in $\mathbb{Z}[\omega]$ is 
\[
n=\omega^a\pi_3^\alpha\pi_1^{\alpha_1}\overline{\pi}_1^{\alpha_1}\cdots
\pi_k^{\alpha_k}\overline{\pi}_k^{\alpha_k}q_1^{\beta_1}\cdots q_l^{\beta_l},
\]
for some $a\in\mathbb{N}$. It follows that if we have the factorization
\begin{equation}\label{nfact}
n=(A+B\omega)(A+B\overline{\omega})=A^2+AB+B^2,
\end{equation}
then the multiplicity of $\pi_i$ in $A+B\omega$ can be any integer $m$ with $0\leq m \leq\alpha_i$, and the multiplicity of $\overline{\pi}_i$ must then be $\alpha_i-m$.\\
A necessary and sufficient condition for $n$ to be representable in the form \eqref{nfact} is that all the $\beta_i$ are even numbers. However, the $q_i$ have no effect on the number of solutions nor on the arguments of the solutions $A+B\omega$.\\
The power $\alpha$ of $3$ in $n$, increases the argument of each solution by a multiple of $\frac{\pi}{6}$, but has no effect on the number of solutions, since the prime factors of $3$ are associated.
Finally, $A+B\omega$ can be multiplied by any unit. This yields that the number of solutions to \eqref{nfact} is  
\begin{equation}\label{rq}
r_Q(n)=6\cdot\prod_{p_i\equiv 1\pod{3}}(\alpha_i+1).
\end{equation}
\begin{example}
The circle in Figure \ref{circle21} has radius $r=21$. The lattice points on the circle therefore correspond to factorizations of
\[
21^2=3^2\cdot7^2
\]
in $\mathbb{Z}[\omega]$. Here
\[
3=(i\sqrt{3})(-i\sqrt{3}), \qquad 7=(3-\omega)(3-\overline{\omega}).
\]
If $\mu=A+B\omega$ lies on the circle, so that $21^2=\mu\overline{\mu}$, then $\mu$ is one of the three numbers
\[
3(3-\omega)^2, \quad 21 \quad \text{or} \quad 3(3-\overline{\omega})^2
\]
multiplied by any of the six units, giving $18$ lattice points.
\end{example}
We have answered the first question on page \pageref{q1}, and consequently devote the rest of the paper to the second one. 
%-----------------------------------------------------------------------
\subsection{Equidistribution through Exponential Sums}
%-----------------------------------------------------------------------
We want to study the distribution of hexagonal lattice points on circles. What we would like to conclude is that the arguments of all the lattice points of a given circle is in some sense evenly distibuted in the interval $[0,2\pi)$. Let us define equidistribution of a sequence in the following intuitive way:
\begin{defn}\label{eqdef}
A sequence $\{x_n\}_{n\in\mathbb{Z}_+}\subseteq [A,B)$ is said to be \emph{equidistributed} if
\[
\lim_{N\to\infty}\frac{\#\{n\leq N;x_n\in[a,b)\}}{N}=\frac{b-a}{B-A}
\]
for all $A\leq a<b\leq B$.
\end{defn}
Equidistribution can be formulated in terms of cancellation in exponential sums as follows:
\begin{thm}[Weyl's Criterion]
The sequence $\{x_n\}$ is equidistributed in $[0,2\pi)$ if and only if
\[
\lim_{N\to\infty} \frac1N \sum_{n=1}^N e^{-ikx_n} = 0
\]
for all integers $k\neq 0$.
\begin{proof}
The proof is based upon the observation that a characteristic function of an interval can be approximated by trigonometric polynomials. (See \cite{kuipers}, Ch.~1, Th.~2.1.)
\end{proof}
\end{thm}
To what extent are the arguments of hexagonal lattice points on circles equidistributed? In Section \ref{sectbad} we will show that there exist arbitrarily large circles with arbitrarily poorly distributed lattice points, so we do not have equidistribution for all circles. We can however prove that lattice points on circles are equidistributed \emph{on average}. Inspired by Weyl's Criterion, we introduce exponential sums:
\begin{defn}\label{defsnk} 
Let $n\in\mathbb{N}$, $A\in\mathbb{Z}\setminus\{0\}$. Define
\[
S(n,A)=\sum_{\substack{\mu\in\mathfrak{O}\\|\mu|^2=n}}e^{iA\arg \mu}.
\]
\end{defn}
Then the main result of this paper, states that the arguments of hexagonal lattice points on circles are equidistributed on average in the following sense:
%-----------------------------------------------------------------------
\begin{thm}\label{main}
If $6 \nmid A$ then $S(n,A)=0$.
If $A\neq 0$ and $6\mid A$ we have for every $\delta >1-\frac2\pi$ , as $x\to\infty$ and $A=O_\delta(e^{\sqrt{\log x}})$,
\[
\frac 1x\sum_{n\le x}|S(n,A)|\ll_\delta (\log x)^{-\delta}.
\]
\end{thm} 
%-----------------------------------------------------------------------
One can interpret this result in the following way: the trivial estimate for $\sum_{n\leq x}|S(n,A)|$ is
\[
\sum_{n\leq x}|S(n,A)|\leq
\sum_{n\leq x}r_Q(n) \sim \pi x.
\]
By Theorem \ref{main}, however, the sum is $o(x)$, ergo the terms are in some sense smaller than expected. This is enough for the application to the Boltzmann Equation. \\
We will also show that the arguments of the hexagonal \emph{primes}, or equivalently, \emph{prime ideals}, are equidistributed. We can define a unique argument of every prime ideal in the following way:     
%-----------------------------------------------------------------------
\begin{defn}\label{thetap}
Let $\mathfrak{p}$ be a prime ideal. If $\mathfrak{p}=(\alpha)$ there is a unique $\alpha'$ among the associates of $\alpha$ satisfying $-\frac\pi 6\leq\arg{\alpha'}\leq\frac\pi 6$. We then define $\theta_\mathfrak{p}=\arg \alpha'$.
\end{defn}
%-----------------------------------------------------------------------
We will then prove
%-----------------------------------------------------------------------
\begin{thm}\label{ku8'''}
The sequence $\{\theta_\mathfrak{p}\}$, where $\mathfrak{p}$ ranges over all prime ideals of $\mathfrak{O}$, ordered by growing $N(\mathfrak{p})$, is equidistributed in $[-\frac{\pi}{6},\frac{\pi}{6})$.
\end{thm}
%-----------------------------------------------------------------------
\subsection{Outline of what to come}
In Section \ref{sectCh} we will introduce two number theoretic concepts related to the exponential sums described above - Hecke characters and Hecke L-functions. The definition of the Hecke characters is in general much more complicated than in this special case, where matters are simplified by the unique factorization in the ring of integers of $K$. The main result of Section \ref{sectCh} is the functional equation for the Hecke L-functions, a special case of the ones derived by Hecke for general algebraic number fields in \cite{hecke2}.

In Section \ref{sectAs} we will use analytic methods to estimate a sum of Hecke characters over prime ideals. The methods are essentially those that are used in the classical proofs of the Prime Number Theorem and the Prime Number Theorem for Arithmetic Progressions (with error terms). In his paper \cite{kubilius}, Kubilius proved these results for the Gaussian number field $\mathbb{Q}(i)$.

In Section \ref{sectEq} we will prove equidistribution on average of lattice points on circles, using a lemma on mean values of multiplicative functions. The method is similar to the one used in the papers by K\'atai-K\"ornyei \cite{kornyei} and Erd\H os-Hall \cite{erdos} concerning the distribution of Gaussian integers.

Section \ref{sectbad} contains a theorem about the density of prime ideals in circle sectors, and a construction of circles with many points that are poorly distributed.

Finally, in Section \ref{sectDis} we give an estimate for the discrepancy, which is another measure of equidistribution. 
%-----------------------------------------------------------------------
\section{Characters and L-functions}\label{sectCh}
\subsection{Definitions and Basic Properties}
Let $K=\mathbb{Q}(\sqrt{-3})$ and let $\mathfrak{O}=\mathbb{Z}[\omega]$, the ring of integers in K. If $\mathfrak{a}$ is an ideal in $\mathfrak{O}$, we will write 
\[
\mathfrak{a}\unlhd\mathfrak{O}.
\]
A \emph{Hecke character} (modulo $(1)$) of the number field $K$ is a function $\chi^{6a}$ from the ideals of $\mathfrak{O}$ to the unit circle defined by:
$$
\chi^{6a}((\mu))=\chi^{6a}(\mu)=\left(\frac{\mu}{|\mu|}\right)^{6a}, a\in\mathbb{Z}.
$$  
This is well defined since every ideal in $\mathfrak{O}$ is principal, and since $(\mu)=(\nu)$ implies $\mu=\varepsilon\nu$, where $\varepsilon^{6}=1$. \\
The \emph{Hecke L-function} $L(s,\chi^{6a})$ is defined by
$$
L(s,\chi^{6a})=\sum_{\mathfrak{a}\unlhd\mathfrak{O}}\frac{\chi^{6a}(\mathfrak{a})}{N(\mathfrak{a})^{s}}=\sum_{(\mu)}\frac{\mu^{6a}}{|\mu|^{2s+6a}}=\frac{1}{6}\sum_{\mu\in\mathfrak{O}}\frac{\mu^{6a}}{|\mu|^{2s+6a}},\qquad \re(s)>1.
$$
(Following the conventional notation, we use the complex variable $s=\sigma+it$.)
%-----------------------------------------------------------------
\begin{rmk}
When $a=0$, the L-function reduces to the Dedekind zeta function 
$$
\zeta_{K}(s)=\sum_{\mathfrak{a}\unlhd\mathfrak{O}}\frac{1}{N(\mathfrak{a})^{s}}.
$$
\end{rmk}
%-----------------------------------------------------------------
We must of course prove that the series defining the L-functions converge. First of all we note that the series
$$
\sum_{\mathfrak{p}}\frac{1}{N(\mathfrak{p})^{\sigma}},
$$
where $\mathfrak{p}$ runs through all prime ideals in $\mathfrak{O}$, converges for $\sigma>1$. Indeed:
$$
\sum_{\mathfrak{p}}\frac{1}{N(\mathfrak{p})^{\sigma}}\leq2\sum_{p}p^{-\sigma}<2\sum_{n}n^{-\sigma}<\infty.
$$
Now we prove
%-----------------------------------------------------------------
\begin{prop} 
The series
$$
\sum_{\mathfrak{a}\unlhd\mathfrak{O}}\frac{\chi^{6a}(\mathfrak{a})}{N(\mathfrak{a})^{s}}
$$
converges absolutely for $\sigma>1$ and uniformly for $\sigma\geq1+\delta>1$, making $L(s,\chi^{6a})$ an analytic function for $\sigma>1$.
Furthermore, for $\sigma>1$, we have the following Euler Product representation:
\begin{equation}\label{ep}
L(s,\chi^{6a})=\prod_{\mathfrak{p}\unlhd\mathfrak{O}}\left(1-\frac{\chi^{6a}(\mathfrak{p})}{N(\mathfrak{p})^{s}}\right)^{-1},
\end{equation}
where the infinite product is absolutely convergent.
\begin{proof}
The general factor of the right side of \eqref{ep} is
$$
\left(1-\frac{\chi^{6a}(\mathfrak{p})}{N(\mathfrak{p})^{s}}\right)^{-1}=\exp\left(\sum_{m=1}^{\infty}\frac{\chi^{6a}(\mathfrak{p}^{m})}{mN(\mathfrak{p})^{ms}}\right).
$$
This is bounded by
$$
\exp\left(\sum_{m=1}^{\infty}\frac{1}{N(\mathfrak{p})^{m\sigma}}\right)=\exp\left(\frac{N(\mathfrak{p})^{-\sigma}}{1-N(\mathfrak{p})^{-\sigma}}\right)\leq\exp(2N(\mathfrak{p})^{-\sigma})
$$
since $N(\mathfrak{p})\geq2$ for all prime ideals $\mathfrak{p}$. Thus, by the above remark the product in \eqref{ep} is absolutely convergent for $\sigma>1$.

Now,
\begin{equation}\label{starsum}
\begin{split}
\prod_{N(\mathfrak{p})\leq x}\left(1-\frac{\chi^{6a}(\mathfrak{p})}{N(\mathfrak{p})^{s}}\right)^{-1}&=\prod_{N(\mathfrak{p})\leq x}\left(1+\frac{\chi^{6a}(\mathfrak{p})}{N(\mathfrak{p})^{s}}+\frac{\chi^{6a}(\mathfrak{p})^{2}}{N(\mathfrak{p})^{2s}}+\ldots\right)\\
&=\sum_{\mathfrak{a}}^{\star}\frac{\chi^{6a}(\mathfrak{a})}{N(\mathfrak{a})^{s}}\\
&=\sum_{N(\mathfrak{a})\leq x}\frac{\chi^{6a}(\mathfrak{a})}{N(\mathfrak{a})^{s}}+\sum_{N(\mathfrak{a})>x}^{\star}\frac{\chi^{6a}(\mathfrak{a})}{N(\mathfrak{a})^{s}},
\end{split}
\end{equation}
where the star indicates that $\mathfrak{a}$ runs through those ideals whose divisors all have norm $\leq x$. In particular, letting $a=0$ and $s=\sigma>1$ in \eqref{starsum}, we see that
$$
\sum_{N(\mathfrak{a})\leq x}\frac{1}{N(\mathfrak{a})^{\sigma}}<\prod_{\mathfrak{p}}\left(1-\frac{1}{N(\mathfrak{p})^{\sigma}}\right)^{-1},
$$
establishing the convergence of
$$
\sum_{\mathfrak{a}}\frac{1}{N(\mathfrak{a})^{\sigma}}.
$$
Also, by \eqref{starsum}
\begin{equation*}
\left|\prod_{N(\mathfrak{p})\leq x}\left(1-\frac{\chi^{6a}(\mathfrak{p})}{N(\mathfrak{p})^{s}}\right)^{-1}-\sum_{N(\mathfrak{a})\leq x}\frac{\chi^{6a}(\mathfrak{a})}{N(\mathfrak{a})^{s}}\right|\leq\sum_{N(\mathfrak{a})>x}^{\star}\frac{1}{N(\mathfrak{a})^{\sigma}}
\to0, 
\end{equation*}
when $x\to\infty$, which completes the proof. 
\end{proof}
\end{prop} 
%-----------------------------------------------------------------
\begin{prop}\label{ld}
For $\re(s)>1$
$$
-\frac{L'(s,\chi^{6a})}{L(s,\chi^{6a})}=\sum_{\mathfrak{p}}\sum_{m=1}^{\infty}\frac{\chi^{6a}(\mathfrak{p}^{m})\log N(\mathfrak{p})}{N(\mathfrak{p})^{ms}}
$$
\begin{proof}
By \eqref{ep} we have, for $\sigma>1$
$$
-\log L(s,\chi^{6a})=\sum_{\mathfrak{p}}\log\left(1-\frac{\chi^{6a}(\mathfrak{p})}{N(\mathfrak{p})^{s}}\right)=-\sum_{\mathfrak{p}}\sum_{m=1}^{\infty}\frac{\chi^{6a}(\mathfrak{p}^{m})}{mN(\mathfrak{p})^{ms}}.
$$
Differentiation yields
$$
-\frac{L'(s,\chi^{6a})}{L(s,\chi^{6a})}=\sum_{\mathfrak{p}}\sum_{m=1}^{\infty}\frac{\chi^{6a}(\mathfrak{p}^{m})\log N(\mathfrak{p})}{N(\mathfrak{p})^{ms}}.
$$
\end{proof}
\end{prop}
%-----------------------------------------------------------------
\begin{rmk}
From now on we assume that $a\neq0$.
\end{rmk}
%-----------------------------------------------------------------
\subsection{A Theta Formula}
The analytic continuation of the L-function to the whole plane will be proved through a transformation formula for a so called theta-function. 
We will give the theta formula that was proved by Hecke \cite{hecke2}, in the particular case of the field $K=\mathbb{Q}(\sqrt{-3})$. The idea behind the proof is to use two-dimensional Poisson Summation. We shall recall some Fourier theory. For variables in $\mathbb{R}^2$ it is always to be understood that $x=(x_1,x_2)$, $y=(y_1,y_2)$ and so on. 
%-----------------------------------------------------------------
\begin{defn}
A function $f\in C^\infty(\mathbb{R}^2)$ is called a \emph{Schwartz function} if it approaches zero faster than any inverse power of $x$ as $|x|\to\infty$, as do all of its derivatives. 
\end{defn}
%-----------------------------------------------------------------
\begin{defn}
We let $\langle \cdot,\cdot \rangle$ be the standard scalar product in $\mathbb{R}^2$:
\[
\langle x,y \rangle:=x_1y_1+x_2y_2.
\]
\end{defn}
%-----------------------------------------------------------------
\begin{defn}
For a Schwartz function $f$ we define the \emph{Fourier transform} $\hat{f}$ by \[
\hat{f}(y)=\int_{\mathbb{R}^2}f(x)e^{2\pi i\langle x,y \rangle}dx.
\]
 \end{defn}
%-----------------------------------------------------------------
\begin{thm}[Poisson Summation Formula]
Let $f$ be a Schwartz function. Then for every $x\in \mathbb{R}^2$ we have
\[
\sum_{m\in\mathbb{Z}^2}f(x+m)=\sum_{m\in\mathbb{Z}^2}\hat{f}(m)e^{2\pi i\langle m,x \rangle}.
\]
\begin{proof}
See Lang (\cite{lang}, XIII,\S1)
\end{proof}
\end{thm}
%-----------------------------------------------------------------
The Fourier transform of the function
\[
h(x)=e^{-\pi \langle x,x \rangle}
\]
is particularly simple - it is easily seen that we have
\begin{equation}\label{h}
\hat{h}=h.
\end{equation}
%-----------------------------------------------------------------
\begin{defn}
Let $f$ be a Schwartz function, $B$ a non-singular real matrix. Define
\[
f_B(x)=f(Bx).
\]
\end{defn}
%-----------------------------------------------------------------
Obviously, $f_B$ is then also a Schwartz function. The Fourier transform is easily found:
%-----------------------------------------------------------------
\begin{lemma}\label{fB}
We have
\[
\hat{f}_B(y)=\frac{1}{|B|}\hat{f}((B^{-1})^T y),
\]
where $|B|$ is the absolute value of the determinant of $B$.
\begin{proof}
We have
\[
\hat{f}_B(y)=\int f(Bx)e^{-2\pi i\langle x,y \rangle}dx.
\]
By the change of variables $z=Bx$ we get
\begin{align*}
\hat{f}_B(y)&=\frac{1}{|B|}\int f(z)e^{-2\pi i\langle B^{-1}z,y \rangle}dz\\
&=\frac{1}{|B|}\int f(z)e^{-2\pi i\langle z,(B^{-1})^T y \rangle}dz\\
&=\frac{1}{|B|}\hat{f}((B^{-1})^T y).
\end{align*}
\end{proof}
\end{lemma}  
%-----------------------------------------------------------------
\begin{defn}
Let $Q(x)$ be a positive definite quadratic form given by a matrix $A$ (that is, $Q(x)=\langle Ax,x \rangle$). We define the Schwartz function $f_Q(x)$ by
\[
f_Q(x)=e^{-\pi Q(x)}.
\]
\end{defn}
%-----------------------------------------------------------------
\begin{lemma}\label{fQ}
We have
\[
\hat{f}_Q(y)=\frac{1}{\sqrt{|A|}}f_{Q'}(y),
\]
where $Q'(x)=\langle A^{-1}x,x \rangle$.
\begin{proof}
Since $A$ is positive definite, there exists a real matrix $B$ such that $B^2=A$. Then $Q(x)=\langle Ax,x \rangle = \langle Bx,Bx \rangle$, so that we have
\[
f_Q(x)=h_B(x),
\]
where $h(x)=e^{-\pi \langle x,x \rangle}$. Thus, by Lemma \ref{fB} and \eqref{h},
\[
\hat{f}_Q(y)=\hat{h}_B(y)=\frac{1}{|B|}\hat{h}(B^{-1}y)=\frac{1}{\sqrt{|A|}}h(B^{-1}y)=\frac{1}{\sqrt{|A|}}e^{-\pi Q'(y)}.
\]
\end{proof}
\end{lemma}
%-----------------------------------------------------------------
Now we prove a formula from which the desired theta formula will follow.
%-----------------------------------------------------------------
\begin{thm}\label{theta1}
Let $x_1,x_2$ be real variables and define
\begin{equation}\label{u1u2}
\begin{cases}
u_1=x_1+\omega x_2,\\
u_2=x_1+\overline{\omega} x_2.
\end{cases}
\end{equation}
Moreover, let $t$ be strictly positive. Then we have the following formula:
\begin{equation}\label{theta1eq1}
\sum_{\mu\in\mathfrak{O}}e^{-2\pi t(\mu+u_1)(\mu+u_2)}
=\frac{1}{\sqrt{3}t}\sum_{\nu\in\mathfrak{O}}e^{-\frac{2\pi}{3t}|\nu|^2+\frac{2\pi}{\sqrt{3}}(\overline{\nu}u_2-\nu u_1)}.
\end{equation}
\begin{rmk}
For real $x_1,x_2$ we have $\overline{u}_1=u_2$. Later we will extend this result to complex $x_1,x_2$ in which case this relation does not hold in general.
\end{rmk}
\begin{proof}
Define the positive quadratic form $Q(y)$ by
\[
Q(y)=2t|y_1+y_2\omega|^2=2t(y_1^2+y_1y_2+y_2^2).
\]
Then $Q(y)=\langle Ay,y \rangle$, where
\[
A=2t
\begin{pmatrix}
1&\frac12\\
\frac12&1
\end{pmatrix}.
\]
Since
\[
A^{-1}=\frac{2}{3t}
\begin{pmatrix}
1&-\frac12\\
-\frac12&1
\end{pmatrix},
\]
we get
\[
Q'(y)=\langle A^{-1}y,y \rangle =\frac{2}{3t}(y_1^2-y_1y_2+y_2^2).
\]
Now let us define
\[
F(x)=\sum_{\mu\in\mathfrak{O}}e^{-2\pi t(\mu+u_1)(\overline{\mu}+u_2)}
=\sum_{\mu\in\mathfrak{O}}e^{-2\pi t|\mu+u_1|^2}.
\]
Then we have
\[
F(x)=\sum_{m\in\mathbb{Z}^2}e^{-2\pi t|m_1+m_2\omega+x_1+x_2\omega|^2}
=\sum_{m\in\mathbb{Z}^2}e^{-\pi Q(x+m)}=\sum_{m\in\mathbb{Z}^2}f_Q(x+m).
\] 
Using Poisson Summation and Lemma \ref{fQ}, we get
\[
F(x)=\sum_{m\in\mathbb{Z}^2}\hat{f}_Q(m)e^{2\pi i \langle m,x \rangle}
=\frac{1}{\sqrt{|A|}}\sum_{m\in\mathbb{Z}^2}e^{-\pi Q'(m)+2\pi i \langle m,x \rangle}.
\]
Solving \eqref{u1u2} for $x_1,x_2$ yields
\[
\begin{cases}
x_1=\frac{i}{\sqrt{3}}(\overline{\omega}u_1-\omega u_2),\\
x_2=\frac{i}{\sqrt{3}}(-u_1+u_2).
\end{cases}
\]
Thus we get
\[
\langle m,x \rangle = m_1x_1+m_2x_2
=\frac{i}{\sqrt{3}}\big(u_1(m_1\overline\omega-m_2)-u_2(m_1\omega-m_2)\big).
\]
Put $\nu=m_1\overline\omega-m_2$. Then if $m$ runs through $\mathbb{Z}^2$, $\nu$ runs through $\mathfrak{O}$. Moreover
\[
Q'(m)=\frac{2}{3t}(m_1^2-m_1m_2+m_2^2)=\frac{2}{3t}|\nu|^2
\]
and
\[
\langle m,x \rangle = \frac{i}{\sqrt{3}}(\nu u_1-\overline\nu u_2).
\]
We conclude that
\[
F(x)=\frac{1}{\sqrt{3}t}\sum_{\nu\in\mathfrak{O}}e^{-\frac{2\pi}{3t}|\nu|^2+\frac{2\pi}{\sqrt{3}}(\overline\nu u_2-\nu u_1)},
\]
as stated.
\end{proof}
\end{thm}
%-----------------------------------------------------------------
We will now extend this result to complex numbers $x_1,x_2$. As before we let
\[
\begin{cases}
u_1=x_1+\omega x_2,\\
u_2=x_1+\overline{\omega} x_2.
\end{cases}
\]
Note that when we allow arbitrary complex values for $x_1$ and $x_2$, $u_1$ and $u_2$ become independent complex variables.
%-----------------------------------------------------------------
\begin{prop}\label{fanal}
The series
\[
F(x)=\sum_{\mu\in\mathfrak{O}}e^{-2\pi t(\mu+u_1)(\mu+u_2)} \quad \text{and} \quad G(x)=\sum_{\nu\in\mathfrak{O}}e^{-\frac{2\pi}{3t}|\nu|^2+\frac{2\pi}{\sqrt{3}}(\overline{\nu}u_2-\nu u_1)} 
\] 
are absolutely convergent for all $(x_1,x_2)\in\mathbb{C}^2$. Moreover, they converge uniformly in the region $\Omega_R=\{(x_1,x_2)\in\mathbb{C}^2\ :\ \max(|x_1|,|x_2|)<R\}$ for every $R>0$.
\begin{proof}
Let 
\[
\mathcal{P}_\mu(x)=2\pi t(\mu+u_1)(\mu+u_2),
\]
so that
\[
F(x)=\sum_{\mu\in\mathfrak{O}}e^{-\mathcal{P}_\mu(x)}.
\]
Expanding
\[
\mathcal{P}_\mu(x)=2\pi t(|\mu|^2+\mu u_2+\overline\mu u_1+u_1u_2),
\]
we see that, in $\Omega_R$,
$$
\re(\mathcal{P}_\mu(x))>2\pi t(|\mu|^2-4|\mu|R-4R^2) \gg |\mu|^2.
$$
Thus, for all but a finite number of $\mu$,
\[
|e^{-\mathcal{P}_\mu(x)}| \leq e^{-c|\mu|^2},
\]
and $\sum_{\mu}e^{-c|\mu|^2}$ clearly converges, so the uniform convergence of the series $F(x)$ follows. Moreover,
\[
G(x)=\sum_{\nu\in\mathfrak{O}}e^{-\mathcal{Q}_\nu(x)},
\]
where
\[
\mathcal{Q}_\nu(x)=\frac{2\pi}{3t}|\nu|^2-\frac{2\pi}{\sqrt{3}}(\overline{\nu}u_2-\nu u_1).
\]
In $\Omega_R$ we have
\[
\re(\mathcal{Q}_\nu(x))>\frac{2\pi}{3t}|\nu|^2-\frac{4\pi R}{\sqrt{3}}|\nu| \gg |\nu|^2,
\]
so the uniform convergence of $G(x)$ follows by an analogous argument.
\end{proof}
\end{prop}
%-----------------------------------------------------------------------
Thus both sides of \eqref{theta1eq1} define entire functions of the complex variables $x_1,x_2$. Since they agree for real $x_1,x_2$, they must be equal. Thus we have proven
\begin{thm}\label{Theta2}
For $t>0$ and arbitrary complex numbers $u_1,u_2$ we have
\begin{equation}\label{theta2}
\sum_{\mu\in\mathfrak{O}}e^{-2\pi t(\mu +u_1)(\overline{\mu}+u_2)}=\frac{1}{\sqrt{3}t}\sum_{\nu\in\mathfrak{O}}e^{-\frac{2\pi}{3t}|\nu|^2+\frac{2\pi}{\sqrt{3}}(\overline{\nu}u_2-\nu u_1)}.
\end{equation}
\end{thm} 
From Theorem \ref{Theta2} we derive our theta formula:
%-----------------------------------------------------------------------
\begin{defn}
Let
$$
\theta(t,a):=\sum_{\mu\in\mathfrak{O}}\mu^{6a}e^{-\frac{2\pi}{\sqrt{3}}t|\mu|^{2}}
$$
\end{defn}
%-----------------------------------------------------------------------
\begin{thm}\label{thetaformula}
We have
$$
\theta(t,a)=t^{-1-6a}\theta(\frac{1}{t},a)
$$
\begin{proof}
We let $\rho\in\mathfrak{O}$ be arbitrary. Furthermore we introduce the variable $z$, and set
\begin{align*}
u_{1}&=\rho & u_{2}&=z+\overline{\rho}
\end{align*}
Now \eqref{theta2} is equivalent to
\begin{multline}\label{2}
\sum_{\mu\in\mathfrak{O}}\exp\Big\{-2\pi t\big(|\mu+\rho|^{2}+z(\mu+\rho)\big)\Big\}
\\
=\frac{1}{\sqrt{3}t}
\sum_{\nu\in\mathfrak{O}}\exp\Big\{-\frac{2\pi}{3t}|\nu|^{2}+\frac{2\pi}{\sqrt{3}}(\overline{\nu}(z+\overline{\rho})-\nu\rho)\Big\}.
\end{multline}
We differentiate (\ref{2}) $6a$ times with respect to $z$. This yields
\begin{multline}
\sum_{\mu\in\mathfrak{O}}(\mu+\rho)^{6a}\exp\Big\{-2\pi t|\mu+\rho|^{2}+z(\mu+\rho)\Big\}\\
=\left(\frac{1}{\sqrt{3}t}\right)^{1+6a}\sum_{\nu\in\mathfrak{O}}\overline{\nu}^{6a}\exp\Big\{-\frac{2\pi}{3t}|\nu|^{2}+\frac{2\pi}{\sqrt{3}}(\overline\nu(z+\overline{\rho})-\nu\rho)\Big\}. 
\end{multline}
Setting $z=0$, we get
\begin{multline}\label{3}
\sum_{\mu\in\mathfrak{O}}(\mu+\rho)^{6a}\exp\Big\{-2\pi t|\mu+\rho|^{2}\Big\}\\
=\left(\frac{1}{\sqrt{3}t}\right)^{1+6a}\sum_{\nu\in\mathfrak{O}}\overline{\nu}^{6a}\exp\Big\{-\frac{2\pi}{3t}|\nu|^{2}+\frac{2\pi}{\sqrt{3}}(\overline{\nu\rho}-\nu\rho)\Big\}. 
\end{multline}
Put $t=\frac{\tau}{\sqrt{3}}$. As $\mu$ ranges over $\mathfrak{O}$, so does $\mu+\rho$, since $\rho\in\mathfrak{O}$. Thus the left side of (\ref{3}) equals
$$
\sum_{\mu\in\mathfrak{O}}\mu^{6a}e^{-\frac{2\pi}{\sqrt{3}}\tau|\mu|^{2}}=\theta(\tau,a).
$$
We note that
\[
\frac{1}{\sqrt{3}}(\overline{\nu\rho}-\nu\rho)
=-i\frac{2}{\sqrt{3}}\im(\nu\rho)
\]
But $\im(\nu\rho)=m\frac{\sqrt{3}}2$ for some $m\in\mathbb{Z}$, and thus we can neglect the term $\frac{2\pi}{\sqrt{3}}(\overline{\mu\rho}-\mu\rho)$ in the exponent of the right side of \eqref{3}.
Therefore the right side equals
\[
\left(\frac{1}{\tau}\right)^{1+6a}\sum_{\nu\in\mathfrak{O}}\overline{\nu}^{6a}e^{-\frac{2\pi}{\sqrt{3}\tau}|\nu|^2}=
\tau^{-1-6a}\sum_{\nu\in\mathfrak{O}}\nu^{6a}e^{-\frac{2\pi}{\sqrt{3}}\frac{1}{\tau}|\nu|^{2}}
=\tau^{-1-6a}\theta(\frac{1}{\tau},a)
\]
since when $\nu$ runs through $\mathfrak{O}$, so does $\overline{\nu}$. This finishes the proof.
\end{proof} 
\end{thm}
%-------------------------------------------------------------------------
\subsection{The Functional Equation}
We will now see how the L-functions can be extended to the whole complex plane.
\begin{defn}
Let
$$
\xi(s,\chi^{6a})=\left(\frac{\sqrt{3}}{2\pi}\right)^{s}\Gamma(s+3|a|)L(s,\chi^{6a})
$$
\end{defn}
%-------------------------------------------------------------------------
\begin{thm}
$\xi(s,\chi^{6a})$ is entire and satisfies the functional equation
$$
\xi(s,\chi^{6a})=\xi(1-s,\chi^{6a})
$$
\begin{proof}
Since $L(s,\chi^{6a})=L(s,\chi^{6(-a)})$ we can assume that $a$ is positive. For $\mu\in K$ and $\sigma>1$ we have
\[
\Gamma(s+3a)|\mu|^{-2(s+3a)}=\int_{0}^{\infty}e^{-t|\mu|^2}t^{s+3a-1}dt,
\]
and hence
\begin{align*}
\Gamma(s+3a)\frac{\chi^{6a}(\mu)}{N(\mu)^{s}}&=\int_{0}^{\infty}\mu^{6a}e^{-t|\mu|^2}t^{s+3a-1}dt\\
&=\left(\frac{2\pi}{\sqrt{3}}\right)^{s+3a}\int_{0}^{\infty}\mu^{6a}e^{-\frac{2\pi}{\sqrt{3}}v|\mu|^{2}}v^{s+3a-1}dv.
\end{align*}
Thus
\begin{equation}\label{5}
\begin{split}
\Gamma(s+3a)\left(\frac{\sqrt{3}}{2\pi}\right)^{s+3a}L(s,\chi^{6a})&=
\sum_{(\mu)\unlhd\mathfrak{O}}\int_{0}^{\infty}\mu^{6a}e^{-\frac{2\pi}{\sqrt{3}}v|\mu|^2}v^{s+3a-1}dv\\
&=\frac{1}{6}\sum_{\mu\in\mathfrak{O}}\int_{0}^{\infty}\mu^{6a}e^{-\frac{2\pi}{\sqrt{3}}v|\mu|^2}v^{s+3a-1}dv\\
\end{split}
\end{equation}
For $\re(s)$ large enough, \eqref{5} is absolutely convergent, and hence we can change the order of integration and summation to get
\begin{equation}\label{5'}
\begin{split}
\Gamma(s+3a)\left(\frac{\sqrt{3}}{2\pi}\right)^{s+3a}L(s,\chi^{6a})
&=\frac{1}{6}\int_{0}^{\infty}\sum_{\mu\in\mathfrak{O}}\mu^{6a}e^{-\frac{2\pi}{\sqrt{3}}v|\mu|^2}v^{s+3a-1}dv\\
&=\frac{1}{6}\int_{0}^{\infty}\theta(v,a)v^{s+3a-1}dv.
\end{split}
\end{equation}
The right side of \eqref{5'} is
\[
=\frac{1}{6}\int_0^1\theta(v,a)v^{s+3a-1}dv + \frac{1}{6}\int_1^\infty\theta(v,a)v^{s+3a-1}dv.
\]
Using Theorem \ref{thetaformula}, this is
\begin{align*}
&=\frac{1}{6}\int_{0}^{1}\theta(\frac{1}{v},a)v^{s-3a-2}dv + \frac{1}{6}\int_{1}^{\infty}\theta(v,a)v^{s+3a-1}dv\\
&=\frac{1}{6}\int_{1}^{\infty}\theta(v,a)v^{-s+3a}dv + \frac{1}{6}\int_{1}^{\infty}\theta(v,a)v^{s+3a-1}dv.
\end{align*}
Thus we have proven (for $\re(s)$ large enough)
\begin{equation}\label{6}
\xi(s,\chi^{6a})={1\over6}\left(\frac{\sqrt{3}}{2\pi}\right)^{-3a}\int_{1}^{\infty}\theta(v,a)(v^{s+3a-1}+v^{-s+3a})dv.
\end{equation}
But since this integral converges absolutely for all $s$, \eqref{6} represents an analytic continuation of $\xi(s,\chi^{6a})$ to the whole plane. The functional equation also follows, since the right side of \eqref{6} remains unchanged when we replace $s$ by $1-s$.
\end{proof}
\end{thm}
From this theorem we see that $L(s,\chi^{6a})$ is an entire function. Using our knowledge of the Gamma function, we also deduce that $L(s,\chi^{6a})$ must have zeros at $s=-3|a|,-1-3|a|,-2-3|a|,\ldots$\ . These are called the trivial zeros. Moreover, for $\sigma>1$ it can be seen from the Euler product that $L(s,\chi^{6a})\neq0$.\\
There is however an infinite number of non-trivial zeros in the so called critical strip, $0\leq\sigma\leq1$. It is generally believed that all of them lie on the line $\re(s)=\frac12$. This statement is part of the generalized Riemann Hypothesis, and if we assume it, the estimates in Section \ref{sectAs} can be made much sharper. In the next section, however, we will unconditionally narrow down the region in which the non-trivial zeros can appear.
%-----------------------------------------------------------------------
\begin{cor}\label{c4}
In the strip $-\frac 12\leq\sigma\leq 4$
$$
L(s,\chi^{6a})=k_1e^{k_2t}
$$
for some positive constants $k_1,k_2$ (depending on $a$).
\begin{proof}
By \eqref{6},
$$
\Gamma(s+3|a|)\left(\frac{\sqrt{3}}{2\pi}\right)^{s+3a}L(s,\chi^{6a})\ll \int_{1}^{\infty}e^{-{2\pi\over\sqrt{3}}u}u^{3+3a}du=O_a(1),
$$
and hence
\[
L(s,\chi^{6a})\underset{a}{\ll}\left(\frac{2\pi}{\sqrt{3}}\right)^{s+3a}\frac{1}{\Gamma(s+3|a|)}=O_a\left(\frac{1}{|\Gamma(s+3|a|)|}\right).
\]
Stirling's Formula states that, in the angular region $-\pi+\delta<\arg s<\pi+\delta$ for any fixed $\delta>0$, we have as $|s|\to\infty$
\[
\log \Gamma(s)=(s-\frac{1}{2})\log s-s+\frac 12\log 2\pi +O(|s|^{-1}).
\]
and hence
\[
\log \frac 1{\Gamma(s)}=(\frac 12-s)\log s+s+O(1).
\]
In the strip $-\frac 12\leq\sigma\leq 4$, since $-\frac\pi 2<\arg(s+3|a|)<\frac\pi 2$, we get
\begin{align*}
\re\left(\log \frac 1{\Gamma(s+3|a|)}\right)
&=\frac 12\log\big|s+3|a|\big|-(\sigma+3|a|)\log\big|s+3|a|\big|\\
&+t\arg(s+3|a|)+\sigma+3|a|+O(1)\\
&<k_2|t|.  
\end{align*}
The corollary follows.
\end{proof}
\end{cor}

%-----------------------------------------------------------------------
\section{An Asymptotic Formula for $\sum_{N(\mathfrak{p})\le x}\chi^{6a}(\mathfrak{p})$}\label{sectAs}
%-----------------------------------------------------------------------
Following Kubilius \cite{kubilius} and Landau \cite{landau2} we will in this section obtain an estimate for $\sum_{N(\mathfrak{p})\le x}\chi^{6a}(\mathfrak{p})$ in terms of $x$ and $a$, by finding zero-free regions for the Hecke L-functions $L(s,\chi^{6a})$ and estimates for the logarithmic derivatives of $L(s,\chi^{6a})$ in these regions.  For the convenience of the reader we recall the following two lemmas from Landau's book: 
\begin{lemma}[Landau \cite{landau1}, Satz 374]\label{la374}
Let $r>0$. Suppose $f(s)$ is analytic for
$|s-s_{0}|\leq r$. Furthermore, suppose
$$
\left|{f(s)\over f(s_{0})}\right|<e^{M} \qquad \textrm{for } |s-s_{0}|\leq r
$$
and
$$
f(s)\neq0 \qquad \textrm{for } |s-s_{0}|\leq r,\quad \re(s)>\re(s_0)
$$
Then the following holds:
\begin{itemize}
\item[1)]
$$
-\re\left(\frac{f'(s_{0})}{f(s_{0})}\right)<\frac{4M}{r}
$$
\item[2)]
If there is a zero $\rho$ on the line between $s_{0}-{r\over2}$ and $s_{0}$ (exclusive), then
$$
-\re\left(\frac{f'(s_{0})}{f(s_{0})}\right)<\frac{4M}{r}-\frac{1}{s_{0}-\rho}.
$$
\end{itemize}
\end{lemma}
%-----------------------------------------------------------------------
\begin{lemma}[Landau \cite{landau1}, Satz 225]\label{la225}
Let $r>0$. Suppose $f(s)$ is analytic for
$|s-s_{0}|\leq r$ and there satisfies
$$
\re(f(s))\leq M
$$
Let $0<\rho<r$. Then for $|s-s_{0}|\leq\rho$
$$
|f'(s)|\leq\frac{2r}{(r-\rho)^{2}}(|M|+|f(s_{0})|).
$$
\end{lemma}
%-----------------------------------------------------------------------
\begin{rmk}
From now on, $c_{1},c_{2},\dots$ will denote suitably chosen positive constants (independent of $a$).
\end{rmk}
%-----------------------------------------------------------------------
\begin{thm}\label{ku1}
In the strip $-{1\over2}\leq\sigma\leq 4$
$$
|L(s,\chi^{6a})|<c_{1}(1+|a|)^{2}(1+|t|)^{2}.
$$
\begin{proof}
From the functional equation we have
$$
\left|L\left(-{1\over2}+it,\chi^{6a}\right)\right|=\left|\left({\sqrt{3}\over2\pi}\right)^{-2-2it}\right|\left|\frac{\Gamma({3\over2}-it+3|a|)}{\Gamma(-{1\over2}+it+3|a|)}\right|
\left|L\left({3\over2}-it,\chi^{6a}\right)\right|.
$$
Since
$$
\left|L\left(\frac{3}{2}-it,\chi^{6a}\right)\right|\leq\sum_{(\mu)}{1\over N(\mu)^{\frac{3}{2}}}=O(1),
$$
we get
$$
\left|L\left(-{1\over2}+it,\chi^{6a}\right)\right|<c_{2}\left|\frac{\Gamma({3\over2}-it+3|a|)}{\Gamma(-{1\over2}+it+3|a|)}\right|.
$$
Applying twice the functional equation of the Gamma function we get
$$
\left|L\left(-{1\over2}+it,\chi^{6a}\right)\right|<c_{2}\left|{1\over2}-it+3|a|\right|\left|-{1\over2}-it+3|a|\right|<c_{3}(1+|a|)^{2}(1+|t|)^{2}.
$$
Furthermore
$$
|L(4+it,\chi^{6a})|\leq\sum_{(\mu)}{1\over N(\mu)^{4}}=O(1).
$$
Consider now the function
$$
\Lambda(s)=\frac{L(s,\chi^{6a})}{(1+|a|)^{2}(1+s)^{2}}
$$
$\Lambda(s)$ is holomorphic in the strip $-\frac{1}{2}\leq\sigma\leq4$, and since
$$
|\Lambda(s)|\ll\frac{|L(s,\chi^{6a})|}{(1+|a|)^{2}(1+|t|)^{2}},
$$
$\Lambda(s)$ is bounded on $\sigma=-\frac{1}{2}$ and $\sigma=4$ by the above. Furthermore it is $O(e^{ct})$ in the whole strip by Corollary \ref{c4}. Thus, by Phragm\'en-Lindel\"of's Theorem, $\Lambda(s)$ is bounded in the whole strip, and the theorem follows.
\end{proof}
\end{thm}
%-----------------------------------------------------------------------
\begin{lemma}\label{zetapole}
In the strip $1<\sigma<2$,
$$
\zeta_{K}(s)<\frac{2}{\sigma-1}.
$$
\begin{proof}
\begin{equation*}
\begin{split}
\zeta_{K}(s)&=\prod_{\mathfrak{p}}(1-N(\mathfrak{p})^{-s})^{-1}\\
&=(1-3^{-s})^{-1}\prod_{p\equiv1(3)}(1-p^{-s})^{-2}\prod_{p\equiv2(3)}(1-p^{-2s})^{-1}\\
&=\prod_{p}(1-p^{-s})^{-1}\prod_{p\equiv1(3)}(1-p^{-s})^{-1}\prod_{p\equiv2(3)}(1+p^{-s})^{-1}\\
&=\zeta(s)L(s,\chi),
\end{split}
\end{equation*}
where $\zeta(s)$ is the Riemann zeta function, and $L(s,\chi)$ the Dirichlet L-function for the character
$$
\chi(n)=\left\{ \begin{array}{lll}
-1 & \text{if } n\equiv2\mod3\\
0 & \text{if } n\equiv0\mod3\\
1 & \text{if } n\equiv1\mod3\\
\end{array}\right..
$$
Now
$$
|\zeta(s)|\leq\sum_{n=1}^{\infty}n^{-\sigma}\leq 1+\int_{1}^{\infty}u^{-\sigma}du=1+\frac{1}{\sigma-1}<\frac{2}{\sigma-1}.
$$
Furthermore,
$$
L(s,\chi)=\sum_{n=1}^{\infty}\frac{\chi(n)}{n^{s}}=\int_{1}^{\infty}u^{-s}d\{\sum_{n\leq u}\chi(n)\}=s\int_{1}^{\infty}u^{-s-1}\big(\sum_{n\leq u}\chi(n)\big)du,
$$
so since $\bigl|\sum_{n\leq x}\chi(n)\bigr|\leq 1$ for all $x\geq1$,
$$
|L(s,\chi)|\leq 1
$$
and the lemma follows.
\end{proof}
\end{lemma}
\begin{figure}[htbp]
\psfrag{c6}{$c_6\ $}
\psfrag{-c6}{$-c_6\quad$}
\centering
\includegraphics[scale=.5, trim=0 60 0 0, clip]{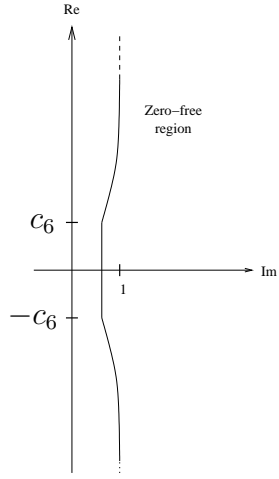}
\caption{A zero-free region for $L(s,\chi^{6a})$}
\end{figure}
%-----------------------------------------------------------------------
The next theorem gives zero-free regions for the Hecke L-functions, of the same form as can be obtained for the Riemann zeta function and the Dirichlet L-functions.
%-----------------------------------------------------------------------
\begin{thm}[Zero-free Region]\label{ku2}
There exist positive constants $c_{5}, c_{6}$ such that $L(s,\chi^{6a})$ has no zeros in the region defined by
$$
\sigma\geq \begin{cases}
1-\frac{1}{c_{5}\log\big((1+|a|)(1+|t|)\big)} & \text{for } |t|\geq c_{6}\\
1-\frac{1}{c_{5}\log\big((1+|a|)(1+|c_{6}|)\big)} & \text{for } |t|\leq c_{6}.
\end{cases}
$$
\end{thm}
\begin{proof}
Logarithmic differentiation of (\ref{ep}) yields
$$
-\frac{L'(s,\chi^{6a})}{L(s,\chi^{6a})}=\sum_{\mathfrak{p}}\sum_{m=1}^{\infty}\frac{\chi^{6ma}(\mathfrak{p})\log N(\mathfrak{p})}{N(\mathfrak{p})^{ms}}
$$
for $\sigma>1$. We recall the definition of $\theta_\mathfrak{p}$ on Page \pageref{thetap} as a unique angle in $[-\frac{\pi}{6},\frac{\pi}{6})$. In terms of $\theta_\mathfrak{p}$ we have, for every $\mathfrak{p}$ and $m$,
\[
\frac{\chi^{6ma}\log N(\mathfrak{p})}{N(\mathfrak{p})^{ms}}
=\frac{\log N(\mathfrak{p})}{N(\mathfrak{p})^{m\sigma}}\exp\big(i(6ma\theta_\mathfrak{p}-mt\log N(\mathfrak{p}))\big).
\]
Using the inequality
$$
3+4\cos\varphi+\cos2\varphi = 2(1+\cos\varphi)^2 \geq 0,
$$
we deduce that
\begin{multline}\label{7}
-3\frac{\zeta'_{K}(\sigma)}{\zeta_{K}(\sigma)}-4\re\left(\frac{L'(\sigma+it,\chi^{6a})}{L(\sigma+it,\chi^{6a})}\right)-\re\left(\frac{L'(\sigma+2it,\chi^{12a})}{L(\sigma+2it,\chi^{12a})}\right)\\
=\sum_{\mathfrak{p}}\frac{\log N(\mathfrak{p})}{N(\mathfrak{p})^{m\sigma}}\Big(3+4\cos\big(6ma\theta_{\mathfrak{p}}-mt\log N(\mathfrak{p})\big)\\
+\cos\big(12ma\theta_{\mathfrak{p}}-2mt\log N(\mathfrak{p})\big)\Big)
\geq 0.
\end{multline}
Let $s_{0}=\rho+i\tau$, where $1<\rho<2$. ($\rho$ is later to be suitably chosen as a function of $\tau$.) In the disk $|s-s_{0}|\leq{3\over2}$ we have, by Theorem \ref{ku1} and Lemma \ref{zetapole}
\begin{equation}
\begin{split}
\left|\frac{L(s,\chi^{6a})}{L(s_{0},\chi^{6a})}\right|&<c_{1}(1+|a|)^{2}(1+|t|)^{2}\left|L(s_{0},\chi^{6a})^{-1}\right|\\
&=c_{1}(1+|a|)^{2}(1+|t|)^{2}\left|\sum_{(\alpha)}\frac{\mu(\alpha)\chi^{6a}(\alpha)}{N(\alpha)^{s_{0}}}\right|\\
&\leq c_{1}(1+|a|)^{2}(1+|t|)^{2}\zeta_{K}(\rho)\\
&<\frac{c_7}{\rho-1}(1+|a|)^{2}(1+|\tau|)^{2}.
\end{split}
\end{equation}
Suppose now that $\mu+i\tau$ is a zero of $L(s,\chi^{6a})$, where $\rho-{3\over4}<\mu<\rho$. Then, applying 2) of Lemma \ref{la374} with
\[
r=\frac 32, \qquad M=\frac{c_7}{\rho-1}(1+|a|)^{2}(1+|\tau|)^{2}
\]
we get
\begin{multline}\label{9}
-\re\left(\frac{L'(\rho+i\tau,\chi^{6a})}{L(\rho+i\tau,\chi^{6a})}\right)<{16\over3}\log((1+|a|)(1+|\tau|))\\
-{8\over3}\log(\rho-1)-{1\over\rho-\mu}+c_{8}.
\end{multline}
Moreover, by 1) of Lemma \ref{la374}
\begin{multline}\label{10}
-\re\left(\frac{L'(\rho+2i\tau,\chi^{12a})}{L(\rho+2i\tau,\chi^{12a})}\right)<{16\over3}\log((1+|a|)(1+|\tau|))\\
-{8\over3}\log(\rho-1)+c_{8}.
\end{multline}
As regards ${\zeta_{K}'\over\zeta_{K}}$ we note that, because of the simple pole at $1$,
\begin{equation}\label{11}
-{\zeta'_{K}(\rho)\over\zeta_{K}(\rho)}<{1\over\rho-1}+c_{10}.
\end{equation}
Now \eqref{7}, \eqref{9}, \eqref{10} and \eqref{11} imply
\begin{equation}\label{12}
{4\over\rho-\mu}<{3\over\rho-1}+{80\over3}\log\big((1+|a|)(1+|\tau|)\big)-{40\over3}\log(\rho-1)+c_{11}.
\end{equation}
We may choose $c_{6}$ sufficiently large to ensure that for $|\tau|\geq c_{6}$
$$
40\log(100\log2(1+|\tau|))+3c_{11}<20\log2(1+|\tau|)
$$
and thus for all $a\neq0$
\begin{multline}\label{13}
40\log\Big(100\log\big((1+|a|)(1+|\tau|)\big)\Big)+3c_{11}\\
<20\log\big((1+|a|)(1+|\tau|)\big). 
\end{multline}
Now we put
$$
\rho=\begin{cases}
1+\frac{1}{100\log\big((1+|a|)(1+|\tau|)\big)} & \text{for }|\tau|\geq c_{6}\\
1+\frac{1}{100\log\big((1+|a|)(1+|c_{6}|)\big)} & \text{for }|\tau|< c_{6}
\end{cases}.
$$
First suppose $|\tau|\geq c_{6}$. Put 
\[
\mathcal{L}=\log\big((1+|a|)(1+|\tau|)\big)
\] 
for short. Then \eqref{12}, multiplied by $3$, becomes
\[
\frac{12}{\rho-\mu} < 980\mathcal{L}+40\log(100\mathcal{L})+3c_{11}.
\]
Thus, by \eqref{13},
\[
\frac{12}{\rho-\mu} < 1000\mathcal{L},
\]
and hence
\[
\mu<1+\frac{1}{100\mathcal{L}}-\frac{12}{1000\mathcal{L}}=1-\frac{1}{500\mathcal{L}}
\]
for all eventual zeros $\mu+i\tau$.

Next suppose $|\tau|<c_{6}$ and put 
\[
\mathcal{L'}=\log\big((1+|a|)(1+|c_{6}|)\big).
\] 
Then
\begin{equation*}
\begin{split}
\rho-\mu&>\frac{12}{\frac{9}{\rho-1}+80\log\big((1+|a|)(1+|\tau|)\big)-40\log(\rho-1)+3c_{11}}\\
&>\frac{12}{980\mathcal{L'}+40\log(100\mathcal{L'})+3c_{11}}.
\end{split}
\end{equation*}
But by \eqref{13} (for $\tau=c_{6}$)
$$
40\log(100\mathcal{L'})+3c_{11}<20\mathcal{L'},
$$
so
\[
\rho-\mu>\frac{12}{1000\mathcal{L'}},
\]
and hence
\[
\mu<\rho-\frac{12}{1000\mathcal{L}'}=1+\frac{1}{100\mathcal{L}'}-\frac{12}{1000\mathcal{L}'}=1-\frac{1}{500\mathcal{L'}}.
\]
We have proven the theorem.
\end{proof}
%-----------------------------------------------------------------------
\begin{lemma}\label{logbound}
On the line $\sigma=2$, 
\[
\log L(s,\chi^{6a})<c_{16}.
\]
\begin{proof}
If $s=2+it$ we have
\begin{align*}
\frac{1}{|L(s,\chi^{6a})|}
&=\prod_{\mathfrak{p}}\left|1-\frac{\chi^{6a}(\mathfrak{p})}{N(\mathfrak{p})^{s}}\right|
\leq\prod_{\mathfrak{p}}\left(1+N(\mathfrak{p})^{-2}\right)\\
&\leq 1+\sum_{\mathfrak{a}}N(\mathfrak{a})^{-2}
=1+\zeta_K(2).
\end{align*}
Thus
\[
\frac{1}{1+\zeta_K(2)} \leq |L(s,\chi^{6a})| \leq \zeta_K(2),
\]
so the logarithm is bounded.
\end{proof}
\end{lemma}
%-----------------------------------------------------------------------
\begin{thm}\label{ku3}
In the region $\Omega$ defined by
$$
3\geq\sigma\geq\begin{cases}
1-\frac{1}{c_{12}\log\big((1+|a|)(1+|t|)\big)} & \text{for }|t|\geq c_{6}\\
1-\frac{1}{c_{12}\log\big((1+|a|)(1+|c_{6}|)\big)} & \text{for }|t|\leq c_{6}
\end{cases},
$$
where $c_{12}>c_{5}$, we have
$$
\left|\frac{L'(s,\chi^{6a})}{L(s,\chi^{6a})}\right|
\leq c_{13}\log^{3}\Big((1+|a|)\big(1+\max(|t_{0}|,c_{6})\big)\Big).
$$
\begin{proof}
Let $c_{14}>c_{5}$. For every $s_{0}=2+it_{0}$ on the line $\re(s)=2$, let $\mathcal{C}_{s_{0}}$ be the circle with center at $s_{0}$ and radius
$$
r=\begin{cases}
1+\frac{1}{c_{14}\log\big((1+|a|)(1+|t_{0}|)\big)} & \text{for }|t_{0}|\geq c_{6}\\
1+\frac{1}{c_{14}\log\big((1+|a|)(1+|c_{6}|)\big)} & \text{for }|t_{0}|\leq c_{6}
\end{cases}.
$$
We may freely assume that $c_{6}>\frac{1}{2}$, so that
$$
(1+|a|)(1+c_{6})>3
$$
Then we have for $s=\sigma+it$ in $\mathcal{C}_{s_{0}}$:
\begin{equation*}\begin{split}
\log\big((1+|a|)(1+|t|)\big)&\leq\log\big((1+|a|)(3+|t_{0}|)\big)\\
&\leq\log\big((1+|a|)(1+|t_{0}|)\big)+\log3
\end{split}\end{equation*}
Thus we get
\begin{equation}\label{14}
\log\big((1+|a|)(1+|t|)\big)<\begin{cases}
2\log\big((1+|a|)(1+|t_{0}|)\big) & \text{if }|t_{0}|\geq c_{6}\\
2\log\big((1+|a|)(1+|c_{6}|)\big) & \text{if }|t_{0}|\leq c_{6}
\end{cases}.
\end{equation}
Now we want to use Lemma \ref{la225} with $f(s)=\log L(s,\chi^{6a})$, which is analytic in the zero-free region of Theorem \ref{ku2}, and thus analytic in $\mathcal{C}_{s_0}$.\\
Also, by Theorem \ref{ku1}
$$
|L(s,\chi^{6a})|<c_{1}(1+|a|)^{2}(1+|t|)^{2}
$$
in  $\mathcal{C}_{s_{0}}$, so by \eqref{14}
\begin{align*}
\re(f(s))
&=\re(\log L(s,\chi^{6a}))\\
&<4\log\big((1+|a|)(1+|t|)\big)+\log c_1\\
&<4\log\Big((1+|a|)\big(1+\max(|t_{0}|,c_{6})\big)\Big)+c_{15}.
\end{align*}
Set
$$
\rho=\left\{\begin{array}{ll}
1+\frac{1}{c_{12}\log\big((1+|a|)(1+|t_{0}|)\big)} & \text{for }|t_{0}|\geq c_{6}\\
1+\frac{1}{c_{12}\log\big((1+|a|)(1+|c_{6}|)\big)} & \text{for }|t_{0}|\leq c_{6}
\end{array}\right.,
$$
where $c_{12}>c_{14}$. Now Lemma \ref{la225}, with $M=4\log((1+|a|)(1+\max(|t_{0}|,c_{6})))+c_{15}$, and Lemma \ref{logbound} imply that in the disk $|s-s_{0}|\leq\rho$ we have
\begin{align*}
\left|\frac{L'(s,\chi^{6a})}{L(s,\chi^{6a})}\right|
&<\frac{2r}{(r-\rho)^{2}}\bigg(4\log\Big((1+|a|)\big(1+\max(|t_{0}|,c_{6})\big)\Big)+c_{15}+|f(s_0)|\bigg)\\
&<\frac{2r}{(r-\rho)^{2}}\bigg(4\log\Big((1+|a|)\big(1+\max(|t_{0}|,c_{6})\big)\Big)+c_{16}\bigg).
\end{align*}
But
\begin{equation*}\begin{split}
\frac{2r}{(r-\rho)^{2}}&=\frac{2r}{\left(\frac{1}{c_{14}}-\frac{1}{c_{12}}\right)^2}\log^{2}\Big((1+|a|)\big(1+\max(|t_{0}|,c_{6})\big)\Big)\\
&\leq\frac{4}{\left(\frac{1}{c_{14}}-\frac{1}{c_{12}}\right)^{2}}\log^{2}\Big((1+|a|)\big(1+\max(|t_{0}|,c_{6})\big)\Big).
\end{split}\end{equation*}
The statement of the theorem follows.
\end{proof}  
\end{thm}
%-----------------------------------------------------------------------
\begin{defn}
Let
$$
K(s,\chi^{6a})=\sum_{\mathfrak{p}}\frac{\chi^{6a}(\mathfrak{p})\log N(\mathfrak{p})}{N(\mathfrak{p})^{s}},\qquad \sigma>1.
$$
\end{defn}
%-----------------------------------------------------------------------
This series is easily seen to be absolutely convergent for $\sigma>1$, and uniformly convergent for $\sigma\geq1+\delta$, so $K(s,\chi^{6a})$ is analytic for $\sigma>1$. We even have:
%-----------------------------------------------------------------------
\begin{lemma}\label{lemmaK}
$K(s,\chi^{6a})$ is analytic in the region $\Omega$ of Theorem \ref{ku3} and satisfies
$$
|K(s,\chi^{6a})|\leq\left\{\begin{array}{ll}
c_{15}\log^{3}\big((1+|a|)(1+|t|)\big) & \text{if }|t|\geq c_{6}\\
c_{15}\log^{3}\big((1+|a|)(1+|c_{6}|)\big) & \text{if }|t|\leq c_{6}
\end{array}\right.
$$
there.
\begin{proof}
By Proposition \ref{ld} we have, for $\sigma>1$
\begin{equation}\label{K}
K(s,\chi^{6a})=-\frac{L'(s,\chi^{6a})}{L(s,\chi^{6a})}-\sum_{\mathfrak{p}}\sum_{m=2}^{\infty}\frac{\chi^{6a}(\mathfrak{p})\log N(\mathfrak{p})}{N(\mathfrak{p})^{ms}}.
\end{equation}
But the logarithmic derivative is analytic in $\Omega$, and the sum on the right is absolutely convergent for $\sigma>\frac{1}{2}$ and uniformly convergent for $\sigma\geq\frac{1}{2}+\delta,\ \delta>0$, since
\begin{equation*}\begin{split}
\sum_{\mathfrak{p}}\sum_{m=2}^{\infty}\frac{\log N(\mathfrak{p})}{N(\mathfrak{p})^{m(\frac{1}{2}+\delta)}}&=\sum_{\mathfrak{p}}\frac{\log N(\mathfrak{p})}{N(\mathfrak{p})^{\frac{1}{2}+\delta}(N(\mathfrak{p})^{\frac{1}{2}+\delta}-1)}\\
&\leq 2\sum_{p}\frac{\log p^{2}}{p^{\frac{1}{2}+\delta}(p^{\frac{1}{2}+\delta}-1)}\\
&<4\sum_{n=2}^{\infty}\frac{\log n}{n^{\frac{1}{2}+\delta}(n^{\frac{1}{2}+\delta}-1)}<\infty.
\end{split}\end{equation*}
Thus, \eqref{K} constitutes an analytic continuation of $K(s,\chi^{6a})$ to $\Omega$ (since clearly $\Omega$ lies to the right of the line $\sigma=\frac{1}{2}$). From \eqref{K} it is also clear that, in $\Omega$,
$$
\left|K(s,\chi^{6a})+\frac{L'(s,\chi^{6a})}{L(s,\chi^{6a})}\right|<c_{17}.
$$
Thus in $\Omega$ we have
\begin{equation*}\begin{split}
|K(s,\chi^{6a})|&<c_{13}\log^{3}\Big((1+|a|)\big(1+\max(|t_{0}|,c_{6})\big)\Big)+c_{17}\\
&<c_{15}\log^{3}\Big((1+|a|)\big(1+\max(|t_{0}|,c_{6})\big)\Big).
\end{split}\end{equation*}
\end{proof}
\end{lemma}
%-----------------------------------------------------------------------
\begin{lemma}\label{int}
$$
\int_{2-i\infty}^{2+i\infty}\frac{x^{s}}{s^{2}}ds=\left\{\begin{array}{ll}
0 & \text{if }0<x<1\\
2\pi i\log x & \text{if }x\geq1
\end{array}\right.
$$
\end{lemma}
\begin{proof}
Let $x>0$ be fixed. The function $\frac{x^s}{s^2}$ is analytic in the whole plane, except for a double pole in the point $s=0$ with residue $\log x$.
\begin{figure}[htbp]
\psfrag{2+iR}{$2+iR$}
\psfrag{2-iR}{$2-iR$}
\psfrag{gR}{$\gamma_R$}
\psfrag{wR}{$\omega_R$}
\begin{minipage}[t]{0.5\linewidth}
\centering
\includegraphics[scale=.5]{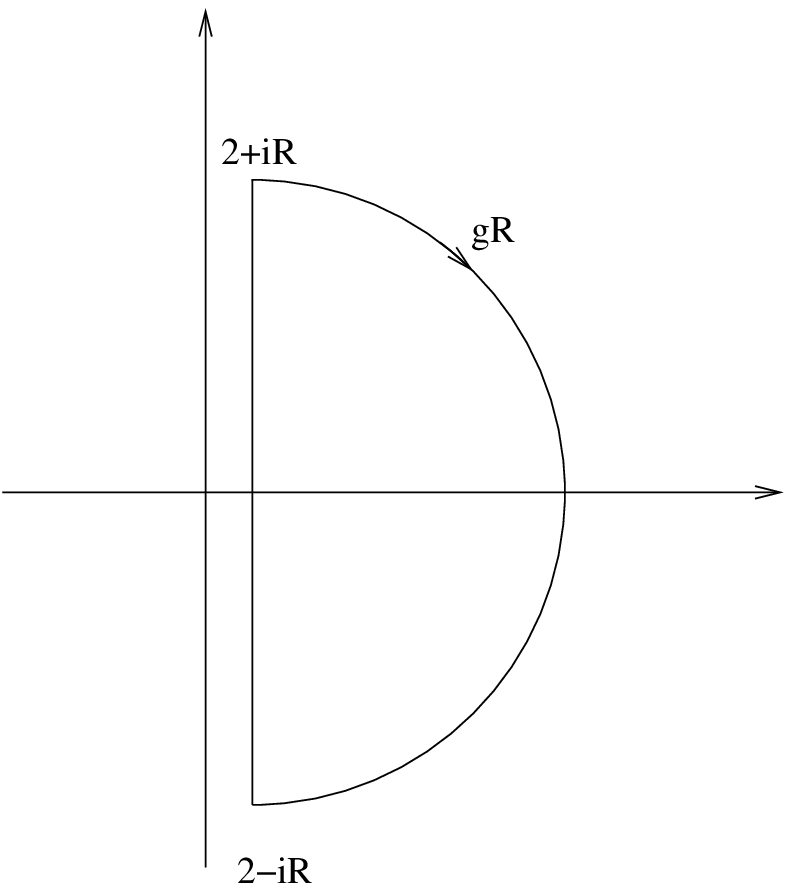}
\caption{}
\label{path1}
\end{minipage}%
\begin{minipage}[t]{0.5\linewidth}
\centering
\includegraphics[scale=.5]{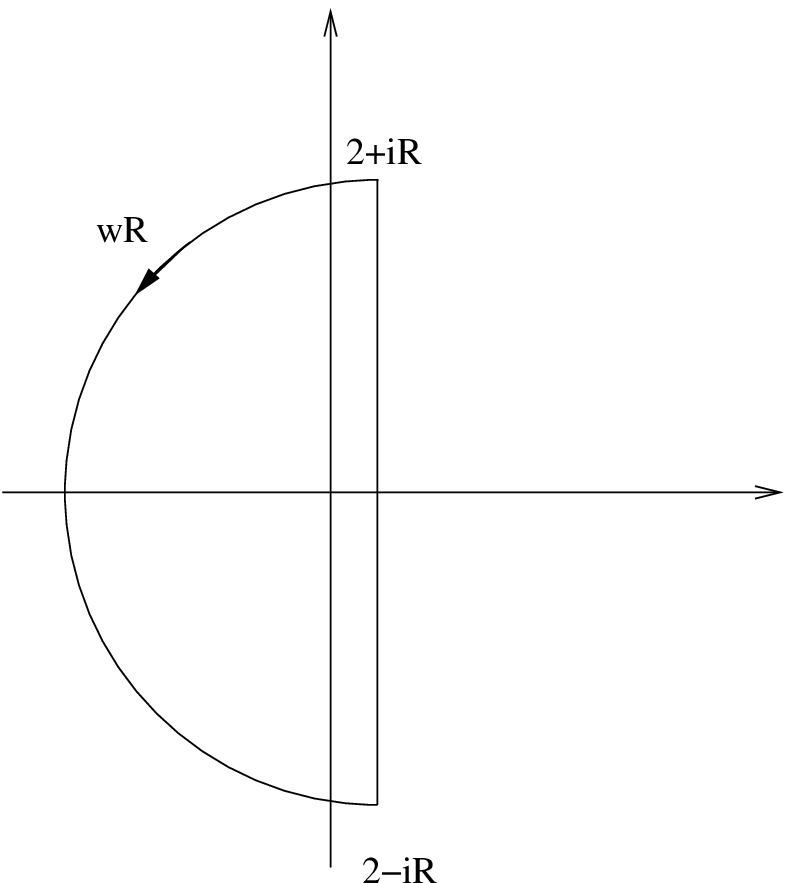}
\caption{}
\label{path2}
\end{minipage}
\end{figure}
Assume first that $0<x<1$. Using the integration contour of Figure \ref{path1} we have by Cauchy's Theorem
\begin{equation*}
\int_{2-iR}^{2+iR}\frac{x^s}{s^2}ds+\int_{\gamma_R}\frac{x^s}{s^2}ds=0.
\end{equation*}
But
\[
\left|\int_{\gamma_R}\frac{x^s}{s^2}ds\right|\leq \pi\frac{x^2}{R^2},
\]
so letting $R\to\infty$ we get
\[
\int_{2-i\infty}^{2+i\infty}\frac{x^s}{s^2}ds=0.
\]
Assume instead that $x\geq1$. We then use the contour of Figure \ref{path2}. Since (if we take $R>2$) the pole lies inside the contour we get by Cauchy's Theorem
\[
\int_{2-iR}^{2+iR}\frac{x^s}{s^2}ds+\int_{\omega_R}\frac{x^s}{s^2}ds=2\pi i\log x.
\]
But
\[
\left|\int_{\omega_R}\frac{x^s}{s^2}ds\right|\leq \pi\frac{x^2}{R^2},
\]
so again we let $R\to\infty$, yielding
\[
\int_{2-i\infty}^{2+i\infty}\frac{x^s}{s^2}ds=2\pi i\log x.
\]
\end{proof}
%-----------------------------------------------------------------------
Now we are ready to estimate sums of characters. We pave the way to $\sum\chi(\mathfrak{p})$ by estimating two weighted sums, with the weights $\log N(\mathfrak{p})\log\frac{x}{N(\mathfrak{p})}$ and $\log N(\mathfrak{p})$, respectively. For the first one, we will use the integral of Lemma \ref{int} to pick out the portion of the L-series for which $N(\mathfrak{p})\leq x$. We will then shift the integration path a bit to the left, close to the zero-free region.
%-----------------------------------------------------------------------
\begin{thm}\label{ku4}
For $x>1$
$$
\sum_{N(\mathfrak{p})\leq x}\chi^{6a}(\mathfrak{p})\log N(\mathfrak{p})\log\frac{x}{N(\mathfrak{p})} \ll xe^{-c_{18}\frac{\log x}{\log(1+|a|)+\sqrt{\log x}}}\log^{3}(1+|a|).
$$
\begin{proof}
By Lemma \ref{int}
\begin{equation}\label{ku41}\begin{split}
\frac{1}{2\pi i}\int_{2-i\infty}^{2+i\infty}\frac{x^{s}}{s^{2}}K(s,\chi^{6a})ds&=\sum_{\mathfrak{p}}\chi^{6a}(\mathfrak{p})\log N(\mathfrak{p})\frac{1}{2\pi i}\int_{2-i\infty}^{2+i\infty}\frac{\left(\frac{x}{N(\mathfrak{p})}\right)^{s}}{s^{2}}ds\\
&=\sum_{N(\mathfrak{p})\leq x}\chi^{6a}(\mathfrak{p})\log N(\mathfrak{p})\log\frac{x}{N(\mathfrak{p})},
\end{split}\end{equation}
which is the sum we want to approximate.
Now let $\omega$ be the curve defined by
$$
\sigma=\left\{\begin{array}{ll}
1-\frac{1}{c_{12}\log\big((1+|a|)(1+|t|)\big)} & \text{for }|t|\geq c_{6}\\
1-\frac{1}{c_{12}\log\big((1+|a|)(1+|c_{6}|)\big)} & \text{for }|t|\leq c_{6}
\end{array}\right.,\qquad -\infty\leq t \leq\infty.
$$
I claim that
\begin{equation}\label{pathshift}
\int_{2-i\infty}^{2+i\infty}\frac{x^{s}}{s^{2}}K(s,\chi^{6a})ds=\int_{\omega}\frac{x^{s}}{s^{2}}K(s,\chi^{6a})ds.
\end{equation}
To see this, consider for large $T$ the contour $\gamma_{T}$ defined in the following manner (see Figure \ref{path3}):
\begin{quote}
from $2-iT$ to $2+iT$ in a straight line,\\
from $2+iT$ to $1-\frac{1}{c_{12}\mathcal{L}}+ iT$ in a straight line,\\
from $1-\frac{1}{c_{12}\mathcal{L}}+ iT$ to $1-\frac{1}{c_{12}\mathcal{L}}- iT$ along $\omega$\\
and from $1-\frac{1}{c_{12}\mathcal{L}}- iT$ to $2-iT$ in a straight line. 
\end{quote}
(Here we have written for short $\mathcal{L}=\log\big((1+|a|)(1+|T|)\big)$.)\\
\begin{figure}[htbp]
\psfrag{T}{$T$}
\psfrag{-T}{$-T$}
\psfrag{w}{$\omega$}
\psfrag{gT}{$\gamma_T$}
\psfrag{c6}{$c_6$}
\psfrag{-c6}{$-c_6$}
\centering
\includegraphics[scale=.5]{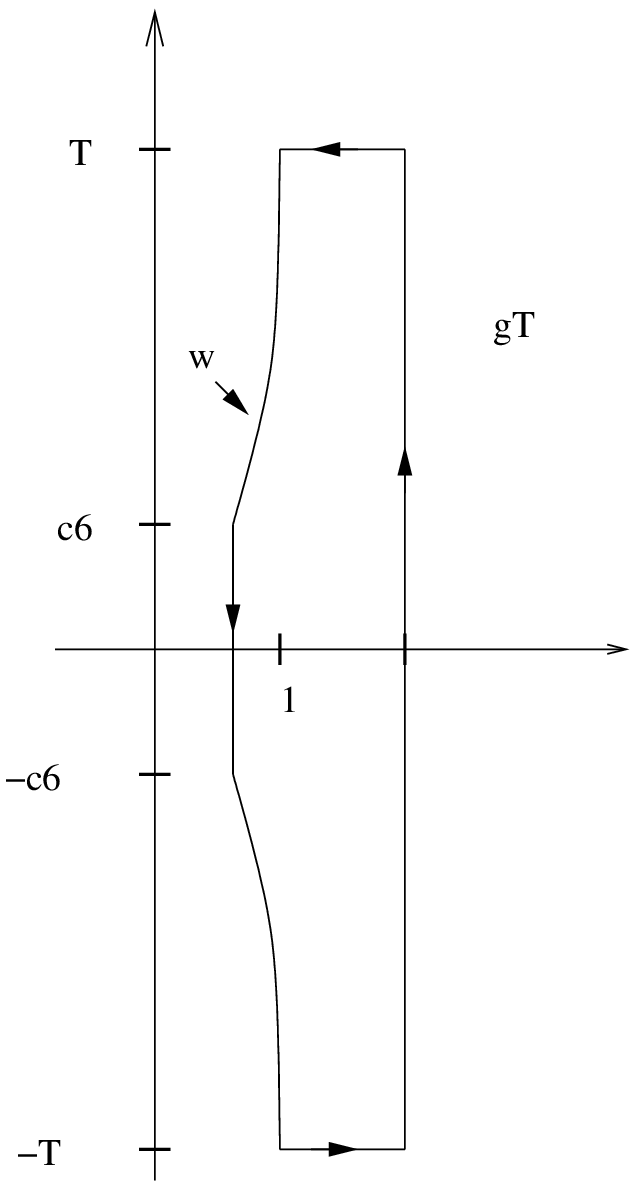}
\caption{}
\label{path3}
\end{figure}
Since the integrand is analytic inside and on $\gamma_{T}$, by Cauchy's Theorem
\begin{equation}\label{cauchy}\begin{split}
0&=\int_{\gamma_{T}}\frac{x^{s}}{s^{2}}K(s,\chi^{6a})ds\\
&=\left(\int_{2-iT}^{2+iT}+\int_{2+iT}^{1-\frac{1}{c_{12}\mathcal{L}}+ iT}+\int_{1-\frac{1}{c_{12}\mathcal{L}}+ iT}^{1-\frac{1}{c_{12}\mathcal{L}}- iT}+\int_{1-\frac{1}{c_{12}\mathcal{L}}- iT}^{2-iT}\right)\frac{x^{s}}{s^{2}}K(s,\chi^{6a})ds.
\end{split}\end{equation}
But by Lemma \ref{lemmaK}
$$
\left|\int_{1-\frac{1}{c_{12}\mathcal{L}}\pm iT}^{2\pm iT}\frac{x^{s}}{s^{2}}K(s,\chi^{6a})ds\right|\leq\frac{x^{2}}{T^{2}}c_{15}\log^{3}\big((1+|a|)(1+T)\big),
$$
so letting $T\to\infty$ in \eqref{cauchy} the horisontal integrals vanish, implying \eqref{pathshift}.
Thus we need to approximate the integral along $\omega$. Now, for an arbitrary $\tau>c_{6}$, we have by Lemma \ref{lemmaK}
\begin{equation*}\begin{split}
\int_{\omega}\frac{x^{s}}{s^{2}}K(s,\chi^{6a})ds
&\ll \int_{0}^{c_{6}} \frac{x^{1-\frac{1}{c_{12}\log((1+|a|)(1+c_{6}))}}}{1+t^{2}}\log^{3}\big((1+|a|)(1+c_{6})\big)dt\\
& \quad +\left(\int_{c_{6}}^{\tau}+\int_{\tau}^{\infty}\right)\frac{x^{1-\frac{1}{c_{12}\log((1+|a|)(1+t))}}}{t^{2}}\log^{3}\big((1+|a|)(1+t)\big)dt.
\end{split}\end{equation*}
The right side is
\begin{equation*}\begin{split}
&\ll x^{1-\frac{1}{c_{12}\log((1+|a|)(1+c_{6}))}}\log^{3}(1+|a|)\\
& \quad +x^{1-\frac{1}{c_{12}\log((1+|a|)(1+\tau))}}\int_{1}^{\infty}\frac{\log^{3}\big((1+|a|)(1+t)\big)}{t^{2}}dt\\
& \quad +x\int_{\tau}^{\infty}\frac{\log^{3}\big((1+|a|)(1+t)\big)}{t^{2}}dt.
\end{split}\end{equation*}
Putting together the first two terms we get
\begin{equation*}\begin{split}
& \ll xe^{-\frac{\log x}{c_{12}\log((1+|a|)(1+\tau))}}\log^{3}(1+|a|)\\ 
& \quad +\frac{x}{\tau}\log^{3}\tau\log^{3}(1+|a|).
\end{split}\end{equation*}
Putting $\tau=e^{\sqrt{\log x}}$, we get
$$
\int_{\omega}\frac{x^{s}}{s^{2}}K(s,\chi^{6a})ds \ll x\log^{3}(1+|a|)\left(e^{-c_{19}\frac{\log x}{\log(1+|a|)+\sqrt{\log x}}}+e^{\frac{3}{2}\log\log x-\sqrt{\log x}}\right).
$$
But for large $x$,
$$
\frac{3}{2}\log\log x-\sqrt{\log x}<-c_{20}\sqrt{\log x}<-c_{20}\frac{\log x}{\log(1+|a|)+\sqrt{\log x}},
$$
so, putting $c_{18}=\min(c_{19},c_{20})$ and recalling \eqref{ku41} and \eqref{pathshift}, we deduce that
$$
\sum_{N(\mathfrak{p})\leq x}\chi^{6a}(\mathfrak{p})\log N(\mathfrak{p})\log\frac{x}{N(\mathfrak{p})} \ll xe^{-c_{18}\frac{\log x}{\log(1+|a|)+\sqrt{\log x}}}\log^{3}(1+|a|).
$$
\end{proof}
\end{thm}
%-----------------------------------------------------------------------
\begin{thm}\label{ku5}
For $x>1$
$$
\sum_{N(\mathfrak{p})\leq x}\chi^{6a}(\mathfrak{p})\log N(\mathfrak{p}) \ll xe^{-c_{21}\frac{\log x}{\log(1+|a|)+\sqrt{\log x}}}\log^{3}(1+|a|).
$$
\begin{proof}
Set for short
$$
\delta=\delta(x)=e^{-\frac{1}{2}c_{18}\frac{\log x}{\log(1+|a|)+\sqrt{\log x}}}.$$
With $x$ replaced by $(1+\delta)x$, Theorem \ref{ku4} gives
\begin{equation}\label{zerothly}\begin{split}
\sum_{N(\mathfrak{p})\leq (1+\delta)x}\chi^{6a}(\mathfrak{p}) & \log N(\mathfrak{p})\log\frac{(1+\delta)x}{N(\mathfrak{p})}\\
&\ll (1+\delta)xe^{-c_{18}\frac{\log (1+\delta)x}{\log(1+|a|)+\sqrt{\log (1+\delta)x}}}\log^{3}(1+|a|)\\
&\ll \delta^{2}x\log^{3}(1+|a|).
\end{split}\end{equation}
We will now split the sum on the left side in two parts. First, again using Theorem \ref{ku4}, we have
\begin{equation}\label{firstly}\begin{split}
\sum_{N(\mathfrak{p})\leq x}\chi^{6a}(\mathfrak{p})\log N(\mathfrak{p})\log\frac{(1+\delta)x}{N(\mathfrak{p})} 
&= \log(1+\delta)\sum_{N(\mathfrak{p})\leq x}\chi^{6a}(\mathfrak{p})\log N(\mathfrak{p})\\
&+ \sum_{N(\mathfrak{p})\leq x}\chi^{6a}(\mathfrak{p})\log N(\mathfrak{p})\log\frac{x}{N(\mathfrak{p})}\\
&= \log(1+\delta)\sum_{N(\mathfrak{p})\leq x}\chi^{6a}(\mathfrak{p})\log N(\mathfrak{p})\\
&+O\big(\delta^{2}x\log^3(1+|a|)\big).
\end{split}\end{equation}
Secondly,
\begin{equation}\label{secondly}\begin{split}
\sum_{x<N(\mathfrak{p})\leq(1+\delta)x}\chi^{6a}(\mathfrak{p})\log N(\mathfrak{p})\log\frac{(1+\delta)x}{N(\mathfrak{p})} 
&\ll \delta x\log((1+\delta)x)\log(1+\delta)\\
&\ll \delta^{2}x\log x,
\end{split}\end{equation}
since the number of terms in this sum is $O(\delta x)$.
By \eqref{firstly} we have
\begin{align*}
\log(1+\delta)\sum_{N(\mathfrak{p})\leq x}\chi^{6a}(\mathfrak{p})\log N(\mathfrak{p})
&=\sum_{N(\mathfrak{p})\leq x}\chi^{6a}(\mathfrak{p})\log N(\mathfrak{p})\log\frac{(1+\delta)x}{N(\mathfrak{p})}\\
&+O\big(\delta^{2}x\log^3(1+|a|)\big)\\
&=\sum_{N(\mathfrak{p})\leq (1+\delta)x}\chi^{6a}(\mathfrak{p})\log N(\mathfrak{p})\log\frac{(1+\delta)x}{N(\mathfrak{p})}\\
&-\sum_{x<N(\mathfrak{p})\leq(1+\delta)x}\chi^{6a}(\mathfrak{p})\log N(\mathfrak{p})\log\frac{(1+\delta)x}{N(\mathfrak{p})}\\
&+O\big(\delta^{2}x\log^3(1+|a|)\big).
\end{align*}
By \eqref{zerothly} and \eqref{secondly}, this is
\begin{align*}
&=O\big(\delta^{2}x\log^3(1+|a|)\big) +O(\delta^2x\log x)+O\big(\delta^{2}x\log^3(1+|a|)\big)\\
&= O\big(\delta^2x\log x\log^3(1+|a|)\big),
\end{align*}
whence
\begin{align*}
\sum_{N(\mathfrak{p})\leq x}\chi^{6a}(\mathfrak{p})\log N(\mathfrak{p}) &\ll \frac{\delta^2}{\log(1+\delta)}x\log x\log^3(1+|a|)\\
&\ll \delta x\log x\log^3(1+|a|)\\
&=xe^{-\frac{1}{2}c_{18}\frac{\log x}{\log(1+|a|)+\sqrt{\log x}}+\log\log x}\log^3(1+|a|)\\
&=xe^{-c_{21}\frac{\log x}{\log(1+|a|)+\sqrt{\log x}}}\log^3(1+|a|).
\end{align*}  
\end{proof}
\end{thm}
%-----------------------------------------------------------------------
The final result of this section now follows easily:
%-----------------------------------------------------------------------
\begin{thm}\label{ku6}
For $x>1$
$$
\sum_{N(\mathfrak{p})\leq x}\chi^{6a}(\mathfrak{p})=\ll xe^{-c_{21}\frac{\log x}{\log(1+|a|)+\sqrt{\log x}}}\log^{3}(1+|a|).
$$
\begin{proof}
Let 
$$
\vartheta(x)=\sum_{N(\mathfrak{p})\leq x}\chi^{6a}(\mathfrak{p})\log N(\mathfrak{p}).
$$
We employ partial summation:
\begin{equation*}\begin{split}
\sum_{N(\mathfrak{p})\leq x}\chi^{6a}(\mathfrak{p})
&=\sum_{2\leq m\leq x}\frac{\vartheta(m)-\vartheta(m-1)}{\log m}\\
&=\sum_{2\leq m\leq x}\vartheta(m)\left(\frac{1}{\log m}-\frac{1}{\log(m+1)}\right)+\frac{\vartheta(x)}{\log([x]+1)}.
\end{split}\end{equation*}
By Theorem \ref{ku5}, this is
\begin{multline*}
\ll \sum_{2\leq m\leq x}me^{-c_{21}\frac{\log m}{\log(1+|a|)+\sqrt{\log m}}}\log^{3}(1+|a|)\left(\frac{1}{\log m}-\frac{1}{\log(m+1)}\right)\\
+ xe^{-c_{21}\frac{\log x}{\log(1+|a|)+\sqrt{\log x}}}\log^{3}(1+|a|).
\end{multline*}
But since $xe^{-c_{21}\frac{\log x}{\log(1+|a|)+\sqrt{\log x}}}$ is monotone increasing for sufficiently large $x$,
\begin{align*}
&\sum_{2\leq m\leq x}me^{-c_{21}\frac{\log m}{\log(1+|a|)+\sqrt{\log m}}}\log^{3}(1+|a|)\left(\frac{1}{\log m}-\frac{1}{\log(m+1)}\right)\\ 
&\ll xe^{-c_{21}\frac{\log x}{\log(1+|a|)+\sqrt{\log x}}}\log^{3}(1+|a|)\sum_{2\leq m\leq x}\left(\frac{1}{\log m}-\frac{1}{\log(m+1)}\right)\\
&\ll xe^{-c_{21}\frac{\log x}{\log(1+|a|)+\sqrt{\log x}}}\log^{3}(1+|a|)\left(\frac{1}{\log 2}-\frac{1}{\log([x]+1)}\right)\\
&\ll xe^{-c_{21}\frac{\log x}{\log(1+|a|)+\sqrt{\log x}}}\log^{3}(1+|a|).
\end{align*}
The theorem follows.
\end{proof}
\end{thm}
%-----------------------------------------------------------------------
\section{Equidistribution of lattice points on circles}\label{sectEq}
%-----------------------------------------------------------------------
\subsection{A Lemma on multiplicative functions}
%-----------------------------------------------------------------------
\begin{defn}A function $f:\mathbb{N}\to\mathbb{C}$ is called \emph{multiplicative} if 
\[
f(mn)=f(m)f(n)
\] 
whenever $(m,n)=1$.
\end{defn}
%-----------------------------------------------------------------------
We will prove a form of the so called \emph{Halberstam-Richert inequality}, giving a bound for $\sum_{n\leq x}f(n)$ via a bound for $\sum_{p\leq x}\frac{f(p)}{p}$. This version appears in K\'atai \cite{katai}. 
%-----------------------------------------------------------------------
\begin{lemma}\label{katai}
Let $f(n)$ be a nonnegative multiplicative function, and assume that 
\[
f(p^\alpha)=O(\alpha)
\]
for every prime power $p^\alpha, \alpha\geq 1$. Then we have
\[
\sum_{n\leq x}f(n) \ll \frac{x}{\log x}\exp\left\{\sum_{p\leq x}\frac{f(p)}{p}\right\}.
\]
\begin{proof}
Let
\begin{align*}
&A(x):=\sum_{n\leq x}f(n), &B(x):=\sum_{n\leq x}f(n)\log n.
\end{align*}
Then we have
\begin{align*}
B(x)&=\sum_{n\leq x}f(n)\sum_{q^\alpha\|n}\log q^\alpha\\
&= \sum_{\substack{q^\alpha h\leq x\\(q^\alpha,h)=1}}f(h)f(q^\alpha)\log q^\alpha,
\end{align*}
since f is multiplicative. Hence
\begin{equation}\label{bx}
B(x)=\sum_{h}f(h)\sum_{q^\alpha\leq\frac{x}{h}}f(q^\alpha)\log q^\alpha.
\end{equation}
But for an arbitrary $y\geq 2$ we get, using the assumption on $f$,
\[
\sum_{q^\alpha\leq y}f(q^\alpha)\log q^\alpha 
\ll \sum_{q^\alpha\leq y}\alpha^2\log q
=\sum_{\alpha}\alpha^2\sum_{q\leq y^{\frac{1}{\alpha}}}\log q.
\]
By the Prime Number Theorem (see for example \cite{davenport}, \S 18)
\[
\sum_{q\leq y^{\frac{1}{\alpha}}}\log q = O(y^{\frac{1}{\alpha}}),
\]
and thus
\[
\sum_{q^\alpha\leq y}f(q^\alpha)\log q^\alpha
\ll \sum_{\alpha}\alpha^2y^{\frac{1}{\alpha}}
= y + \sum_{\alpha>1}\alpha^2y^{\frac{1}{\alpha}}.
\]
For every $\alpha$ in the sum we have $2^\alpha\leq y$, so 
\[
\sum_{\alpha>1}\alpha^2y^{\frac{1}{\alpha}}
\leq \sum_{1<\alpha\leq\frac{\log y}{\log 2}}\alpha^2y^{\frac{1}{\alpha}}
\leq \sum_{1<\alpha\leq\frac{\log y}{\log 2}}\left(\frac{\log y}{\log 2}\right)^2 y^{\frac{1}{2}}
\ll y^{\frac{1}{2}}\log^3y,
\]
since the number of terms in the sum is $O(\log y)$. Thus
\[
\sum_{q^\alpha\leq y}f(q^\alpha)\log q^\alpha
=y+O(y^\frac{1}{2}\log^3y)=O(y).
\]
Inserting this into \eqref{bx} yields
\begin{equation}\label{bx2}
B(x)\ll\sum_{h\leq x}f(h)\frac{x}{h}=x\sum_{h\leq x}\frac{f(h)}{h}.
\end{equation}
But since $f$ is nonnegative and multiplicative, we have
\begin{align*}
\sum_{h\leq x}\frac{f(h)}{h}
&\leq \prod_{p\leq x}\left(1+\frac{f(p)}{p}+\frac{f(p^2)}{p^2}+\ldots\right)\\
&\leq \prod_{p\leq x}\left(1+\frac{f(p)}{p}\right)\left(1+\frac{f(p^2)}{p^2}+\frac{f(p^3)}{p^3}\ldots\right).  
\end{align*}
We claim that
\[ 
\prod_{p}\left(1+\frac{f(p^2)}{p^2}+\frac{f(p^3)}{p^3}\ldots\right)<\infty.
\]
Indeed, this product converges absolutely if and only if the series 
\[
\sum_{p}\left(\frac{f(p^2)}{p^2}+\frac{f(p^3)}{p^3}+\ldots\right)
\]
converges absolutely, and by the assumption on $f$ this series is 
\[
\ll \sum_p\sum_{k=2}^\infty\frac{k}{p^k}<\infty.
\]
Thus
\begin{align*}
\sum_{h\leq x}\frac{f(h)}{h}
&\ll \prod_{p\leq x}\left(1+\frac{f(p)}{p}\right)
= \exp\left\{\sum_{p\leq x}\log\left(1+\frac{f(p)}{p}\right)\right\}\\
&= \exp\left\{\sum_{p\leq x}\left(\frac{f(p)}{p}+O(p^{-2})\right)\right\}
\ll \exp\left\{\sum_{p\leq x}\frac{f(p)}{p}\right\}.
\end{align*}
so, recalling \eqref{bx2}, we get
\begin{equation}\label{bx3}
B(x)\ll x\exp\left\{\sum_{p\leq x}\frac{f(p)}{p}\right\}.
\end{equation}
Now we can make a crude estimate of A(x): by partial Stieltjes integration we have (assume without loss of generality that $x>2$)
\begin{equation}\label{ax}
A(x)=1+\sum_{2\leq n\leq x}f(n)=1+\int_{2-}^x\frac{1}{\log t}dB(t)=1+\frac{B(x)}{\log x}+\int_{2-}^x\frac{1}{t\log^2t}B(t)dt.
\end{equation}
But by Merten's Theorem (see for example \cite{tenenbaum}, Ch.~I.1, Theorem 9), \eqref{bx3} implies
\[
B(x)\ll x\exp\left\{\sum_{p\leq x}\frac{c}{p}\right\} \ll x(\log x)^c,
\]
and hence
\[
\int_{2-}^x\frac{1}{t\log^2t}B(t)dt
\ll \int_2^x (\log t)^{c-2} dt
\ll \begin{cases}
x & \text{ if } c\leq 2,\\
x(\log t)^{c-1} & \text{ if } c>2.
\end{cases}
\]
Putting these estimates into \eqref{ax} yields
\begin{equation}\label{crude}
A(x)\ll x^{1+\varepsilon}
\end{equation}
for an arbitrary $\varepsilon >0$. Now
\begin{align*}
A(x)-A(\sqrt{x})
&= \sum_{\sqrt{x}<n\leq x}f(n)
\leq \frac{2}{\log x}\sum_{\sqrt{x}<n\leq x}f(n)\log n\\
&\leq \frac{2}{\log x}B(x)
\ll \frac{x}{\log x}\exp\left\{\sum_{p\leq x}\frac{f(p)}{p}\right\}.
\end{align*}
But by \eqref{crude}, 
\[
A(\sqrt{x})=O(x^{1/2+\varepsilon}),
\]
so this term is negligible, and  the lemma follows.
\end{proof}
\end{lemma}
%-----------------------------------------------------------------------
\subsection{Proof of Theorem \ref{main}}\label{ssEq}
%-----------------------------------------------------------------------
We recall the definition of the exponential sums $S(n,A)$ on Page \pageref{defsnk}, and define functions $f_A(n)$ in the following way:
\begin{defn}Let
\[
f_A(n):=\frac{|S(n,A)|}{6}.
\]
\end{defn}
These will prove to be multiplicative functions of $n$, thus allowing us to make use of Lemma \ref{katai}. 
\begin{lemma}
$f_A(n)$ is a multiplicative function. 
\begin{proof}
Suppose that $n$ has the prime factorization of \eqref{nfact} on page \pageref{nfact}. If we exclude the trivial case where one or more of the $\beta_i$ are odd, and consequently $f_A(q_i^{\beta_i})=0=f_A(n)$, we get
\[
S(n,A)=\sum_{m=0}^5\sum_{j_1=1}^{\alpha_1}\cdots\sum_{j_k=1}^{\alpha_k}
\omega^{mA}e^{iA\bigl(\alpha\pi/2+(\alpha_1-2j_1)\theta_{p_1}+\ldots+(\alpha_k-2j_k)\theta_{p_k}\bigr)}.
\]
Since
\begin{equation}\label{omegasum}
\sum_{m=0}^5\omega^{mA}=\begin{cases}
6& \text{if } A\equiv 0\pmod{6},\\
0& \text{otherwise},
\end{cases}
\end{equation}
we see that
\begin{align*}
f_A(n)=\frac{|S(n,A)|}{6}
&=\biggl|\sum_{j_1=1}^{\alpha_1}\cdots\sum_{j_k=1}^{\alpha_k}e^{iA\bigl((\alpha_1-2j_1)\theta_{p_1}+\ldots+(\alpha_k-2j_k)\theta_{p_k}\bigr)}\biggr|\\
&=\biggl|\sum_{j_1=1}^{\alpha_1}e^{iA(\alpha_1-2j_1)\theta_{p_1}}\biggr|\cdots
\biggl|\sum_{j_k=1}^{\alpha_k}e^{iA(\alpha_k-2j_k)\theta_{p_k}}\biggr|\\
&=f_A(p_1)\cdots f_A(p_k),
\end{align*}
so $f_A$ is multiplicative.
\end{proof}
\end{lemma}
Furthermore, by \eqref{omegasum}, $f_A=0$ when $A\not\equiv 0 \pmod{6}$, and in this case all results are trivial, so we henceforth assume that $6\mid A$ and substitute $A$ with $6a$, $a\in\mathbb{Z}$.

We examine the values of $f_{6a}$ for prime powers:
\begin{itemize}
\item $f_{6a}(3^\alpha)=1$.
\item If $q\equiv 2\pmod{3}$ we have
\[
f_{6a}(q^\alpha)=\begin{cases}
1 & \text{if $\alpha$ is even,}\\
0 & \text{if $\alpha$ is odd.}
\end{cases}
\]
\item If $p\equiv 1\pmod{3}$ we have
\[
f_{6a}(p^\alpha)=\biggl|\sum_{j=0}^\alpha e^{6ia(\alpha-2j)\theta_p}\biggr|
\le \alpha+1.
\] 
\end{itemize} 
Thus $f_{6a}$ clearly satisfies the hypothesis of Lemma \ref{katai}. Moreover, we have in particular
\begin{equation}\label{f(p)}
f_{6a}(p)=\begin{cases}
2|\cos 6a\theta_p| & \text{if $p\equiv 1\pmod{3}$},\\
0 & \text{if $p\equiv 2\pmod{3}$},\\
1 & \text{if $p=3$}.
\end{cases}
\end{equation}

Now we need to calculate $\sum_{p\le x}\frac{f_{6a}(p)}{p}$. From Fourier analysis we recall Fej\'er's Theorem (\cite{korner}, Th. 1.5), stating that if $g$ is a continuous function on $[0,2\pi]$, then the arithmetic means of the partial sums of the Fourier series of $g$ converge to $g$ uniformly on $[0,2\pi]$. Thus for each $\varepsilon >0$ there exist $k$ and $a_0,\ldots,a_k$ such that for every $x\in [0,2\pi]$
\[
\left||\cos x|-\sum_{m=0}^k a_m\cos mx\right|\le \varepsilon.
\]
There are only cosine-terms in the Fourier series, since the function $|\cos x|$ is even. Moreover
\[
a_0=\frac{1}{2\pi}\int_0^{2\pi}|\cos x|dx=\frac{2}{\pi}.
\]
By \eqref{f(p)} we get
\begin{equation}\label{sumf(p)/p}\begin{split}
\sum_{p\le x}\frac{f_{6a}(p)}{p}
&=\frac 13 + \sum_{\substack{p\le x\\p\equiv 1\pod 3}}\frac{2|\cos 6a\theta_p|}p\\
&\le \frac 13 + \sum_{\substack{p\le x\\p\equiv 1\pod 3}}\left(\sum_{m=0}^k\frac{2a_m\cos 6am\theta_p}p+\frac {2\varepsilon}p\right)\\
&=\frac 13 +2\left(\frac 2\pi+\varepsilon\right) \sum_{\substack{p\le x\\p\equiv 1\pod 3}}\frac 1p \\
&\quad +\sum_{m=1}^k a_m\left(\sum_{\substack{p\le x\\p\equiv 1\pod 3}}\frac{2\cos 6am\theta_p}p\right).
\end{split}\end{equation}
In order to estimate the sum on the right we will reformulate Theorem \ref{ku6} a bit:
%-----------------------------------------------------------------------
\begin{prop}\label{s4prop1}
If $x>1$, $a\neq 0$ and $|am|=O(e^{\sqrt{\log x}})$ we have
\[
\sum_{\substack{p\le x\\p\equiv 1(3)}}2\cos 6am\theta_p \ll xe^{-c_{24}\sqrt{\log x}}.
\]
\begin{proof}
We take a closer look at the sum
\[
\sum_{N(\mathfrak{p})\le x} \chi^{6am}(\mathfrak{p}).
\]
If $\mathfrak{p}$ is a prime ideal, then either $N(\mathfrak{p})=p$, where $p$ is a split prime, or $N(\mathfrak{p})=q^2$, where $q$ is an inert prime, or $N(\mathfrak{p})=3$. In the first case we have $(p)=\mathfrak{p}\overline{\mathfrak{p}}=(\pi_p)(\overline{\pi}_p)$. Thus, if $p\equiv 1\pmod{3}$ we have
\[
\sum_{N(\mathfrak{p})=p}\chi^{6am}(\mathfrak{p})=\chi^{6am}(\pi_p)+\chi^{6am}(\overline{\pi}_p)
=2\cos6am\theta_p.
\]
If $q\equiv 2\pmod{3}$ we have
\[
\sum_{N(\mathfrak{p})=q^2}\chi^{6am}(\mathfrak{p})=\chi^{6am}(q)=1.
\]
Finally
\[
\sum_{N(\mathfrak{p})=3}\chi^{6am}(\mathfrak{p})=\chi^{6am}(\pi_3)=(-1)^{am}.
\]
We conclude that
\begin{equation}\label{s4prop1eq1}
\sum_{N(\mathfrak{p})\le x} \chi^{6am}(\mathfrak{p})
=\sum_{\substack{p\le x\\p\equiv 1(3)}}2\cos 6am\theta_p+
\sum_{\substack{q^2\le x\\q\equiv 2(3)}}1 \ \pm 1.
\end{equation}
Moreover, if $1+|am|\le Ce^{\sqrt{\log x}}$ we have
 \[
 \frac{\log x}{\log(1+|am|)+\sqrt{\log x}} \ge C'\sqrt{\log x},\qquad 0<C'<1,
 \]
 so by Theorem \ref{ku6}
 \begin{align*}
 \sum_{N(\mathfrak{p})\le x}\chi^{6am}(\mathfrak{p})
 &\ll xe^{-c_{21}\frac{\log x}{\log(1+|am|)+\sqrt{\log x}}}\log^3(1+|am|)\\
 &\ll xe^{-c_{22}\sqrt{\log x}}(\log x)^\frac 32\\
 &\le xe^{-c_{23}\sqrt{\log x}}.
 \end{align*}
Thus by \eqref{s4prop1eq1}
\[
\sum_{\substack{p\le x\\p\equiv 1(3)}}2\cos 6am\theta_p 
\ll xe^{-c_{23}\sqrt{\log x}} + x^\frac 12
\ll xe^{-c_{24}\sqrt{\log x}},
\]
and the proposition follows.
\end{proof}
\end{prop}
%-----------------------------------------------------------------------
Note the restriction on $a$. A simple partial integration now yields the estimate we want: 
%-----------------------------------------------------------------------
\begin{prop}
If $x>1$, $a\neq 0$ and $|am|=O(e^{\sqrt{\log x}})$ we have
\begin{equation}\label{s4prop2}
\sum_{\substack{p\le x\\p\equiv 1(3)}}\frac{2\cos 6am\theta_p}p=O(1).
\end{equation}
\begin{proof}
By Proposition \ref{s4prop1}
\begin{align*}
\sum_{\substack{p\le x\\p\equiv 1(3)}}\frac{2\cos 6am\theta_p}p
&=\int_{5-}^x t^{-1} \ d\bigg\{\sum_{\substack{p\le t\\p\equiv 1(3)}}2\cos 6am\theta_p\bigg\}\\
&\ll O(1) + e^{-c_{24}\sqrt{\log x}}+ \int_{5}^x t^{-1}e^{-c_{24}\sqrt{\log t}}dt.
\end{align*}
Now, by the change of variables $u=\sqrt{\log t}$, we have
\begin{align*}
0\leq\int_{5}^x t^{-1}e^{-c_{24}\sqrt{\log t}}dt
&\leq \int_1^\infty t^{-1}e^{-c_{24}\sqrt{\log t}}dt
=\int_0^{\infty} 2ue^{-c_{24}u}du
=\frac{2}{c_{24}^2},
\end{align*}
and thus
\[
\sum_{\substack{p\le x\\p\equiv 1(3)}}\frac{2\cos 6am\theta_p}p
\ll e^{-c_{24}\sqrt{\log x}}+ O(1)=O(1).
\]
\end{proof}
\end{prop}
%-----------------------------------------------------------------------
From the Prime Number Theorem for Arithmetic Progressions (\cite{davenport}, Ch. 20) it easily follows that
\begin{equation}\label{merten2}
\sum_{\substack{p\le x\\ p\equiv 1(3)}}\frac 1 p =\frac 1 2 \log\log x+O(1).
\end{equation}
Inserting \eqref{s4prop2} and \eqref{merten2} into \eqref{sumf(p)/p}, we get
\begin{align*}
\sum_{p\le x}\frac{f_{6a}(p)}{p}
&\le \frac{1}{3}+\left(\frac 2 \pi+\varepsilon\right) \log\log x + O\left(k\max_{1\le m\le k}|a_m|\right)\\
&= \left(\frac 2 \pi+\varepsilon\right) \log\log x + O_\varepsilon(1),
\end{align*}
uniformly in $a$, provided $a\neq0$ and $a=O_\varepsilon(e^{\sqrt{\log x}})$.

Now, using Lemma \ref{katai}, we conclude that
\[
\sum_{n\le x}f_{6a}(n)
\ll \frac{x}{\log x}\exp\left\{\left(\frac{2}{\pi}+\varepsilon\right)\log\log x+O_\varepsilon(1)\right\}
\ll x(\log x)^{\frac 2\pi+\varepsilon -1},
\]
which proves Theorem \ref{main}.\qed
%-----------------------------------------------------------------------
\section{The Distribution of Hexagonal Primes in Sectors. Bad Circles.}\label{sectbad}
Already in 1903 Landau \cite{landau3} proved the so called Prime Ideal Theorem, giving an asymptotic formula for the number of prime ideals with norm $\leq x$ in an arbitrary number field, then with the error term $O(xe^{-(\log x)^{1/13}})$. The version we give here comes from \cite{landau2}. First we define the logarithmic integral:
\begin{defn}
For $x\geq 2$ we define
\[
\Li(x):=\int_2^x\frac{du}{\log u}.
\]
\end{defn}
\begin{thm}[Prime Ideal Theorem]\label{PIT}
For the number of prime ideals in a number field $K$ with norm $\leq x$, we have
\[
\pi_K(x)=\sum_{N(\mathfrak{p})\leq x}1=\Li(x)+O(xe^{-\frac{b}{\sqrt{n}}\sqrt{\log x}}),
\]
where $n$ is the degree of $K$, and $b$ is a positive constant, independent of $K$.
\end{thm}
%-----------------------------------------------------------------------
We will refine Theorem \ref{PIT} in the particular case of the hexagonal number field $K=\mathbb{Q}(\sqrt{-3})$ to a result measuring the number of prime ideals $\mathfrak{p}\unlhd\mathfrak{O}$ such that $N(\mathfrak{p})\leq x$ and such that $\theta_\mathfrak{p}$ lies in a specified interval $[\varphi_1,\varphi_2]$. Again we follow Kubilius \cite{kubilius}.

\begin{figure}
\centering
\includegraphics[scale=0.85]{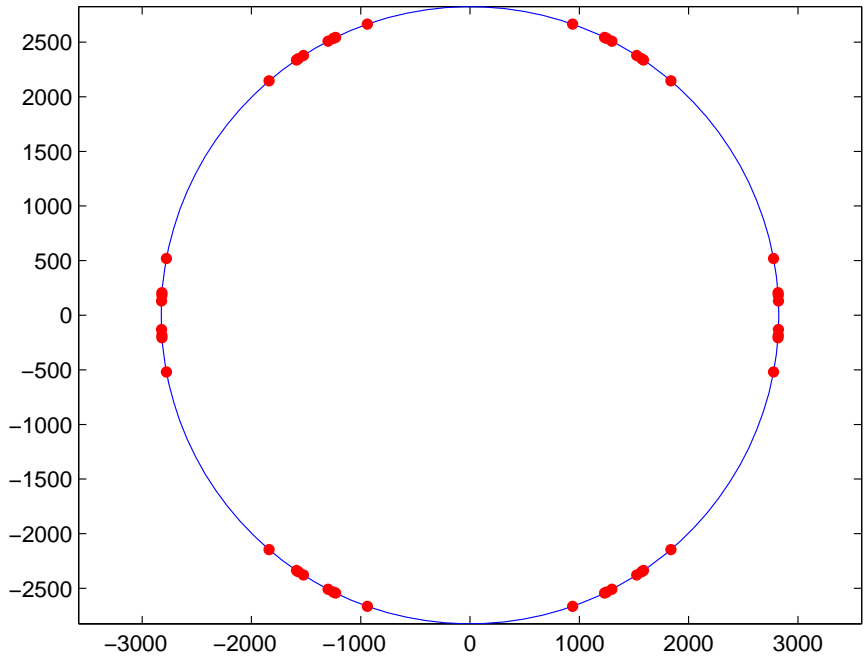}
\caption{Lattice points on the circle with radius $r=\sqrt{7983607}$.}
\label{figbad}
\end{figure}

To this end we will apply the results on the previous section to a Fourier series. We will use a lemma of Vinogradov  to prove that the characteristic function of the interval $[\varphi_1,\varphi_2]$ can be approximated by functions with well-behaved Fourier coefficients.
%-----------------------------------------------------------------------
\begin{lemma}[\cite{vinogradov}, Ch.1, Lemma 12]\label{vino}
Let $r$ be a positive integer, and let $\alpha,\beta,\Delta$ be real numbers satisfying
\[
0<\Delta<1,\quad  \Delta\leq\beta-\alpha\leq 1-\Delta.
\]
Then there exists a periodic function $\psi(x)$, with period $1$, satisfying
\begin{enumerate}
\item $\psi(x)=1$ in the interval $\alpha+\frac\Delta 2\leq x\leq \beta-\frac{\Delta}2$,
\item $\psi(x)=0$ in the interval $\beta+\frac\Delta 2\leq x\leq 1+\alpha-\frac{\Delta}2$,
\item $0\leq\psi(x)\leq 1$ in the remainder of the interval $\alpha-\frac\Delta 2\leq x\leq 1+\alpha-\frac\Delta 2$,
\item $\psi(x)$ has an expansion in Fourier series of the form
\[
\psi(x)=\beta-\alpha+\sum_{m=1}^\infty(a_m \cos 2\pi mx + b_m\sin 2\pi mx),
\]
where
\begin{align*}
&|a_m|,|b_m|\leq 2(\pi m)^{-1},\\
&|a_m|,|b_m|<\frac{2}{\pi m}\left(\frac{r}{\pi m\Delta}\right)^r.
\end{align*}
\end{enumerate}
\end{lemma}
%-----------------------------------------------------------------------
\begin{rmk}
If we instead express the function $\psi(x)$ above in a Fourier series of the form
\[
\psi(x)=\sum_{n=-\infty}^{\infty}c_ne^{2\pi inx},
\]
then we have
\begin{align*}
c_0&=\beta-\alpha, & c_n&=\frac12 (a_n-ib_n), & c_{-n}&=\frac12 (a_n+ib_n), & (n&\geq 1)
\end{align*}
so that 
\begin{align*}
|c_n|&\leq 2(\pi|n|)^{-1},\\
|c_n|&<\frac{2}{\pi|n|}\left(\frac{r}{\pi|n|\Delta}\right)^r.
\end{align*}
\end{rmk}
%-----------------------------------------------------------------------
\begin{lemma}\label{ku7}
Let $\delta >0$ and suppose $2\delta\leq\varphi_2-\varphi_1\leq\frac\pi 3 -2\delta$. Then there exist $\frac{\pi}{3}$-periodic functions $\overline{f}(\varphi)$ and $\underline{f}(\varphi)$ such that
\begin{enumerate}
\item 
\begin{itemize}
\item[] $\overline{f}(\varphi)=1$ if $\varphi_1\leq\varphi\leq\varphi_2$
\item[] $\overline{f}(\varphi)=0$ if $\varphi_2+\delta\leq\varphi\leq\frac\pi 3+\varphi_1-\delta$
\item[] $0\leq\overline{f}(\varphi)\leq 1$ in the rest of the interval $\varphi_1-\delta \leq \varphi \leq \frac\pi 3+\varphi_1-\delta$.
\end{itemize}
\item
\begin{itemize}
\item[] $\underline{f}(\varphi)=1$ if $\varphi_1+\delta\leq\varphi\leq\varphi_2-\delta$
\item[] $\underline{f}(\varphi)=0$ if $\varphi_2\leq\varphi\leq\frac\pi 3+\varphi_1$
\item[] $0\leq\underline{f}(\varphi)\leq 1$ in the rest of the interval $\varphi_1\leq \varphi \leq \frac\pi 3+\varphi_1$.
\end{itemize}
\item if
\[
\overline{f}(\varphi)=\sum_{n=-\infty}^\infty \overline{a}_n e^{6ni\varphi},\qquad \underline{f}(\varphi)=\sum_{n=-\infty}^\infty \underline{a}_n e^{6ni\varphi}
\]
then we have
\begin{align*}
&\overline{a}_0=\frac{3}{\pi}(\varphi_2-\varphi_1+\delta),&\overline{a}_n&\ll \frac{1}{|n|},&\overline{a}_n &\ll \frac{1}{\delta|n|^2}\\
&\underline{a}_0=\frac{3}{\pi}(\varphi_2-\varphi_1-\delta),&\underline{a}_n &\ll \frac{1}{|n|},&\underline{a}_n &\ll \frac{1}{\delta|n|^2}.\\
\end{align*}
\end{enumerate}
\begin{figure}[htbp]
\psfrag{p1-d}{$\scriptstyle \varphi_1-\delta$}
\psfrag{p1}{$\scriptstyle \varphi_1$}
\psfrag{p1+d}{$\scriptstyle \varphi_1+\delta$}
\psfrag{p2-d}{$\scriptstyle \varphi_2-\delta$}
\psfrag{p2}{$\scriptstyle \varphi_2$}
\psfrag{p2+d}{$\scriptstyle \varphi_2+\delta$}
\psfrag{p1+p-d}{$\scriptstyle \varphi_1+\frac{\pi}{3}-\delta$}
\psfrag{p1+p}{$\scriptstyle \varphi_1+\frac{\pi}{3}$}
\psfrag{fo}{$\overline{f}(\varphi)$}
\psfrag{fu}{$\underline{f}(\varphi)$}
\centering
\includegraphics[]{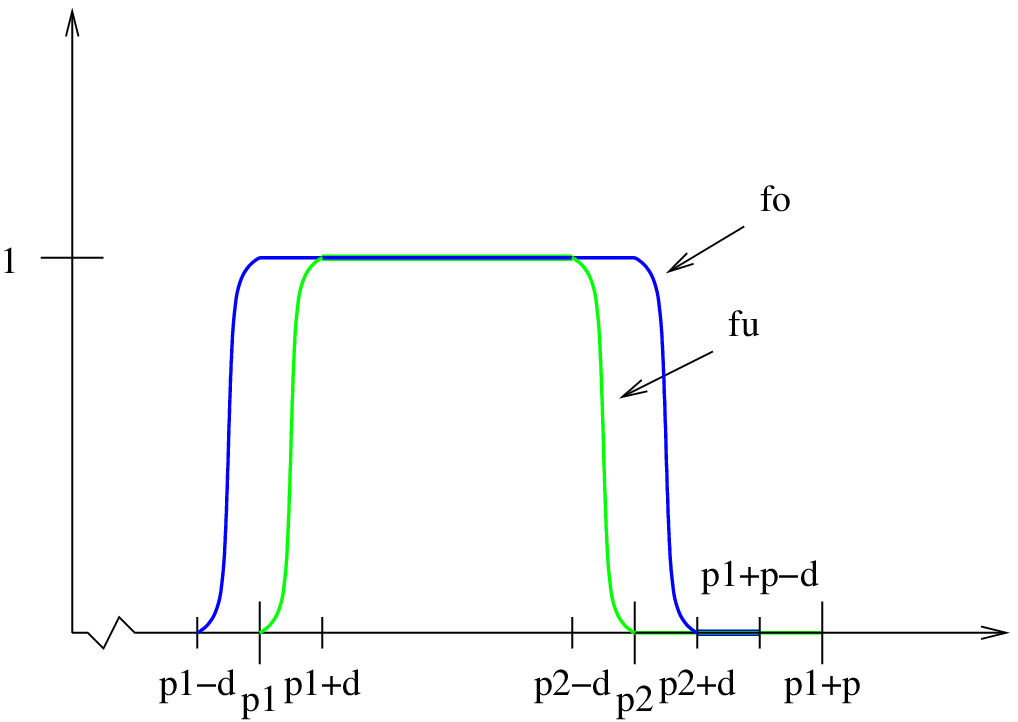}
\caption{}
\end{figure}

\begin{proof}
This follows immediately from Lemma \ref{vino} if we take $x=\frac{3}{\pi}\varphi$ and $r=1$, setting for $\overline{f}$,
\[
\alpha=\frac{3}{\pi}\varphi_1-\frac{3}{2\pi}\delta, \quad \beta=\frac{3}{\pi}\varphi_2+\frac{3}{2\pi}\delta, \quad \Delta=\frac{3}{\pi}\delta,
\]
and for $\underline{f}$
\[
\alpha=\frac{3}{\pi}\varphi_1+\frac{3}{2\pi}\delta, \quad \beta=\frac{3}{\pi}\varphi_2-\frac{3}{2\pi}\delta, \quad \Delta=\frac{3}{\pi}\delta.
\]
\end{proof}
\end{lemma}  
%-----------------------------------------------------------------------
\begin{defn}
Let $-\frac{\pi}{6}\le \varphi_1\le \varphi_2 < \frac{\pi}{6}$. Then we define
\[
\pi_{[\varphi_1,\varphi_2]}(x)=
\sum_{\substack{N(\mathfrak{p})\le x\\
\varphi_1\le \theta_\mathfrak{p} \le \varphi_2}} 1.
\]
\end{defn}
%-----------------------------------------------------------------------
We now prove the theorem about the distribution of prime ideals in circle sectors.
%-----------------------------------------------------------------------
\begin{thm}\label{ku8}
We have
\[
\pi_{[\varphi_1,\varphi_2]}(x)=
\frac3\pi (\varphi_2-\varphi_1) \Li(x)+ O(xe^{-c_{25}\sqrt{\log x}}).
\]
\begin{proof}
Define the functions $\overline{f}(\varphi)$ and $\underline{f}(\varphi)$ as in Lemma \ref{ku7}, with
\[
\delta=e^{-c_{26}\sqrt{\log x}}.
\]
Then we get
\begin{align*}
\pi_{[\varphi_1,\varphi_2]}(x)&=\sum_{\substack{N(\mathfrak{p})\leq x\\
\varphi_1\leq\theta_\mathfrak{p}\leq\varphi_2}}1\leq\sum_{N(\mathfrak{p})\leq x}\overline{f}(\theta_\mathfrak{p})
=\sum_{N(\mathfrak{p})\leq x} \sum_{n=-\infty}^\infty \overline{a}_n\chi^{6n}(\mathfrak{p})\\
&=\quad \overline{a}_0\pi(x) + \sum_{n\neq 0}\overline{a}_n\left(\sum_{N(\mathfrak{p})\leq x}\chi^{6n}(\mathfrak{p})\right).
\end{align*}
Thus, by Theorems \ref{PIT} and \ref{ku6},
\begin{equation}\label{pileq}\begin{split}
\pi_{[\varphi_1,\varphi_2]}(x) 
&\leq \frac3\pi (\varphi_2-\varphi_1+\delta)\left(\Li(x) + O(xe^{-c_{27}\sqrt{\log x}})\right)\\
& \quad + O\left(\sum_{n\neq 1}|\overline{a}_n|xe^{-c_{21}\frac{\log x}{\log(1+|n|)+\sqrt{\log x}}}\log^3(1+|n|)\right). 
\end{split}\end{equation}
Analogously we deduce
\begin{equation}\label{pigeq}\begin{split}
\pi_{[\varphi_1,\varphi_2]}(x)
&\geq \sum_{N(\mathfrak{p})\leq x}\underline{f}(\theta_\mathfrak{p})
= \sum_{N(\mathfrak{p})\leq x} \sum_{n=-\infty}^\infty \underline{a}_n\chi^{6n}(\mathfrak{p})\\
&= \frac3\pi (\varphi_2-\varphi_1-\delta)\left(\Li(x) + O(xe^{-c_{27}\sqrt{\log x}})\right)\\
& \quad + O\left(\sum_{n\neq 1}|\underline{a}_n|xe^{-c_{21}\frac{\log x}{\log(1+|n|)+\sqrt{\log x}}}\log^3(1+|n|)\right).
\end{split}\end{equation}
We examine the sum on the right side of \eqref{pileq}. By the bounds on $|\overline{a}_n|$ in Lemma \ref{ku7} we get, if we split the sum in two parts:
\begin{multline*}
\sum_{n\neq 1}|\overline{a}_n|x e^{-c_{21}\frac{\log x}{\log(1+|n|)+\sqrt{\log x}}}\log^3(1+|n|)\\
\ll x\sum_{1\leq n \leq \delta^{-2}}\frac{\log^3n}{n}e^{-c_{21}\frac{\log x}{\log(1+\delta^{-2})+\sqrt{\log x}}} + x\sum_{n > \delta^{-2}}\frac{\log^3n}{\delta n^2}.
\end{multline*}
For the first part we note that
\[
\log(1+\delta^{-2}) \ll \log \delta^{-2} \ll \sqrt{\log x},
\]
and thus
\begin{align*}
\sum_{1\leq n \leq \delta^{-2}}\frac{\log^3n}{n}e^{-c_{21}\frac{\log x}{\log(1+\delta^{-2})+\sqrt{\log x}}}
&\leq \log^3(\delta^{-2})e^{-c_{28}\sqrt{\log x}} \sum_{1\leq n\leq \delta^{-2}}\frac{1}{n}\\
&\leq e^{-c_{28}\sqrt{\log x}}\log^4(\delta^{-2})\\
&= e^{-c_{28}\sqrt{\log x}}(2c_{26}\sqrt{\log x})^4\\
&\leq e^{-c_{29}\sqrt{\log x}}.
\end{align*}
For the second part we have
\begin{align*}
\sum_{n>\delta^{-2}}\frac{\log^3{n}}{\delta n^2}
&\ll \frac{1}{\delta}\frac{\log^3(\delta^{-2})}{\delta^{-2}}
= \delta\log^3(\delta^{-2})\\
&= e^{-c_{26}\sqrt{\log x}}(2c_{26}\sqrt{\log x})^3
\leq e^{-c_{30}\sqrt{\log x}}.
\end{align*}
Obviously we have the exact same bounds for the corresponding sum in \eqref{pigeq} containing $\underline{a}_n$. Thus \eqref{pileq} and \eqref{pigeq} yield
\[
\left| \pi_{[\varphi_1,\varphi_2]}(x)-\frac3\pi (\varphi_2-\varphi_1)\Li(x) \right|
\ll xe^{-c_{25}\sqrt{\log x}},
\]
and we are done.
\end{proof} 
\end{thm}
%-----------------------------------------------------------------------
Since there is a one-to-one correspondence between prime ideals and hexagonal primes (that is, primes of the number ring $\mathfrak{O}$) in the angular interval $[-\frac{\pi}{6},\frac{\pi}{6})$, we get
%-----------------------------------------------------------------------
\begin{cor}\label{ku8'}
The number of hexagonal primes in the circle sector
\[
\{z;|z|\leq\sqrt{x},\varphi_1\leq\arg z\leq\varphi_2\},\qquad \text{where }
-\frac{\pi}{6}\leq\varphi_1<\varphi_2<\frac{\pi}{6}
\]
is
\[
\frac3\pi (\varphi_2-\varphi_1) \Li(x) + O(xe^{-c_{25}\sqrt{\log x}}).
\]
\end{cor}
%-----------------------------------------------------------------------
We also have
\begin{cor}\label{ku8''}
The same estimate holds if we consider only non-real primes.
\begin{proof}
The number of real hexagonal primes in the sector described above is obviously
\[
\ll \sqrt{x} \ll xe^{-c_{31}\sqrt{\log x}}.
\]
\end{proof}
\end{cor}
%-----------------------------------------------------------------------
Theorem \ref{ku8'''} also follows directly from Theorem \ref{ku8}:
%-----------------------------------------------------------------------
\begin{proof}[Proof of Theorem \ref{ku8'''}]
By Theorem \ref{PIT} and Corollary \ref{ku8'} we have, for $-\frac{\pi}{6}\leq\varphi_1<\varphi_2<\frac{\pi}{6}$
\[
\lim_{x\to\infty}\frac{\pi_{[\varphi_1,\varphi_2]}(x)}{\pi(x)}
=\lim_{x\to\infty}\frac{(1+o(1))\left(\frac3\pi (\varphi_2-\varphi_1) \Li(x)\right)}{(1+o(1))\Li(x)}
=\frac{3}{\pi}(\varphi_2-\varphi_1).
\]
This is what is required in the definition of equidistribution on Page \pageref{eqdef}.
\end{proof}
%-----------------------------------------------------------------------
 We will now show that it is possible to construct arbitrarily ``bad'' circles. We follow Cilleruelo \cite{cilleruelo}.
%-----------------------------------------------------------------------
\begin{thm}\label{cilleruelo}
For every $\varepsilon>0$ and every $k\in\mathbb{N}$ there exists $n\in\mathbb{N}$ such that the circle with radius $\sqrt{n}$ centered at the origin has more than $k$ lattice points, all of which are concentrated on the six arcs
\[
\{\sqrt{n}e^{i(\nu\frac{\pi}{3} + \varphi)};|\varphi|<\varepsilon\},\quad \nu=0,1,\ldots,5.
\] 
\begin{proof}
Let $\varepsilon$ and $k$ be fixed. Choose an integer $m$ such that
\[
m\geq \frac{\log k - \log 6}{\log 2},
\]
and let $0<\delta<\frac{\varepsilon}{m}$. Then by Corollary \ref{ku8'''} we can find $m$ different primes $p_1,\ldots,p_m$ such that $\delta\leq \theta_{p_j} \leq \frac{\varepsilon}{m}$. Set
\[
n=p_1\cdots p_m.
\]
Then the solutions $\alpha\in\mathfrak{O}$ to the equation $n=\alpha\overline{\alpha}$ all have the form
\[
\alpha=\sqrt{n}e^{i(\pm\theta_{p_1}\pm\ldots\pm\theta_{p_m})+ i\nu\frac{\pi}{3}}.
\]
In each case
\[
|\pm\theta_{p_1}\pm\ldots\pm\theta_{p_m}|<\varepsilon,
\] 
and by \eqref{rq} the number of solutions is
\[
r_Q(n)=6\cdot 2^m \geq k,
\]
which proves the result.
\end{proof}
\end{thm}
%-----------------------------------------------------------------------
An example clarifies the method:
%-----------------------------------------------------------------------
\begin{example}
Put
\[
n=7983607=157\cdot 211\cdot 241.
\]
We have
\begin{align*}
157&=12^2+12\cdot1+1^2, & \theta_{157}&\approx0.0692<\frac{\pi}{36},\\
211&=14^2+14\cdot1+1^2, & \theta_{211}&\approx0.0597<\frac{\pi}{36},\\
241&=15^2+15\cdot1+1^2, & \theta_{241}&\approx0.0558<\frac{\pi}{36},
\end{align*}
so that
\[
|\pm\theta_{157}\pm\theta_{211}\pm\theta_{241}|<\frac{\pi}{12}.
\]
Thus the lattice points on the circle with radius $\sqrt{n}$, centered at the origin, all have arguments lying, modulo $\frac{\pi}3$, between $-\frac{\pi}{12}$ and $\frac{\pi}{12}$. Moreover, this circle has $6\cdot 2^3=48$ points. (See Figure \ref{figbad}.)

\end{example}
%-----------------------------------------------------------------------
\section{A Further Measure of Equidistribution}\label{sectDis}
We will in this section give a result that perhaps better justifies the statement of equidistribution on average. In the case of true equidistribution, the ratio of points lying in a specific subinterval would be approximately equal to the ratio between the length of the subinterval and the length of the whole interval. The \emph{discrepancy} $\Delta(n)$ measures how far off this approximation is.  
%-----------------------------------------------------------------------
\begin{defn}
Let
\[
\Delta(n)=\sup_{0\leq \alpha<\beta\leq 2\pi}
\left|\frac{\#\{\zeta\ ;\ |\zeta|^2=n, \arg\zeta\in [\alpha,\beta)\}}{r_Q(n)}-
\frac{\beta-\alpha}{2\pi}\right|.
\]
\end{defn}
%-----------------------------------------------------------------------
We will show that the discrepancy is small for all but a few values of $n$, as $n$ grows large. Obviously we must specify the meaning of ``a few'' - from Section \ref{intro} we know that in some sense only ``a few'' circles with radius $\sqrt{n}$, $n\in\mathbb{N}$, have any lattice points at all (this happens if and only if all prime factors $q\equiv 2\pmod{3}$ occur in even powers in $n$).

We recall that $Q$ is the quadratic form defined by
\[
Q(x,y)=x^2+xy+y^2,
\]
and that $r_Q(n)$ is the number of representations of the integer $n$ by the form $Q$. We put
\[
\mathcal{R}_Q(x):=\{n\leq x\ ;\ r_Q(n)\neq 0\},
\]
\[
B_Q(x):=|\mathcal{R}_Q(x)|.
\]
%-----------------------------------------------------------------------
The asymptotic formula for $B_Q(x)$ was found by Landau \cite{landau4}:
\begin{thm}
There exists a constant $b>0$ such that
\[
B_Q(x)\sim b\frac{x}{\sqrt{\log x}}.
\]
\end{thm}
From now on let $r(n)=r_Q(n)$. Our result will take the form:
%-----------------------------------------------------------------------
\begin{thm}\label{discrepancy}
For almost all $n\in \mathcal{R}_Q(x)$, that is, with the exception of $o(B_Q(x))$ of them, we have
\[
\Delta(n)\leq r_Q(n)^{-\gamma},
\]
if $\gamma<\frac{\log\pi}{\log 2}-1$.
\end{thm}
%-----------------------------------------------------------------------
In proving this we follow K\'atai and K\"ornyei \cite{kornyei}, who gave the analogous result for the square lattice, as did Erd\H os and Hall \cite{erdos}. We will need a result of Erd\H os and Tur\'an \cite{turan} concerning the discrepancy:
\begin{lemma}\label{turanlemma}
Let $\varphi_1,\ldots,\varphi_N\in\mathbb{R}$. Put
\[
Z_k=\frac{1}{N}\sum_{j=1}^N e^{ik\varphi_j}.
\]
Then for an arbitrary $T>0$ we have
\[
\sup_{0\leq\alpha<\beta\leq 2\pi} 
\left|\frac1N\sum_{\alpha\leq\varphi_j<\beta}1-\frac{\beta-\alpha}{2\pi}\right|
\ll \frac1T + \sum_{k=1}^T \frac{|Z_k|}{k}.
\]
\begin{proof}
See \cite{turan}, Th.~III.
\end{proof}
\end{lemma} 
%-----------------------------------------------------------------------
\begin{proof}[Proof of Theorem \ref{discrepancy}]
For $0<\gamma<1$, put
\[
C_\gamma(x)=\sum_{n\leq x}\Delta(n)r(n)^\gamma,
\]
and
\[
M_\gamma(x)=\#\{n\leq x\ ;\ \Delta(n) > r(n)^{-\gamma}\}.
\]
Our aim will be to prove that $C_\gamma(x)=o(B_Q(x))$, since then it follows that $M_\gamma(x)=o(B_Q(x))$. By Lemma \ref{turanlemma} we have for an arbitrary $T>0$
\[
\Delta(n) \ll \frac{1}{T} + \sum_{A=1}^T \frac{\frac{1}{r(n)}|S(n,A)|}{A},
\]
and thus
\begin{align*}
C_\gamma(x) 
&\ll \sum_{n\leq x}\frac{r(n)^\gamma}{T} + \sum_{n\leq x} r(n)^\gamma \sum_{A=1}^T \frac{|S(n,A)|}{A\,r(n)}\\
&\leq \frac{1}{T}\sum_{n\leq x}r(n) + \sum_{A=1}^T \frac{1}{A} \sum_{n\leq x} |S(n,A)|r(n)^{\gamma-1}\\
&\ll \frac{x}{T} + \sum_{A=1}^T \frac{1}{A} \sum_{n\leq x}6^\gamma g_A(n),
\end{align*}
where
\[
g_A(n)=6^{-\gamma}|S(n,A)|r(n)^{\gamma-1} = f_A(n)\left(\frac{r(n)}{6}\right)^{\gamma-1}.
\]
$g_A(n)$ is easily seen to be a multiplicative function of $n$, and since for primes $p$ we have
\[
g_A(p)=2^{\gamma-1}f_A(p), \qquad g_A(p^\alpha)\leq (\alpha+1)^\gamma=O(\alpha),
\]
another application of Lemma \ref{katai} on Page \pageref{katai} yields
\[
\sum_{n\leq x} g_A(n) \ll x(\log x)^{\frac{2^\gamma}{\pi}+\varepsilon-1},
\]
for an arbitrary $\varepsilon>0$. Putting $T=[\log x]+1$ now yields
\[
C_\gamma(x) \ll \frac{x}{\log x} + x(\log x)^{\frac{2^\gamma}{\pi}+\varepsilon-1}\log\log x.
\]  
If we choose $\gamma < \frac{\log\pi}{\log 2}-1$ and $\varepsilon$ sufficiently small, so that $\frac{2^\gamma}{\pi}+\varepsilon <\frac{1}{2}$, we get
\[
C_\gamma(x) = o\left(\frac{x}{\sqrt{\log x}}\right) = o(B_Q(x)),
\]
which completes the proof.
\end{proof}
%-----------------------------------------------------------------------
\section*{Acknowledgements}
I wish to thank my supervisor P\"ar Kurlberg for his enthusiasm and generosity during the work with this thesis. I also express my gratitude to Prof.~Jonas Kubilius for providing a copy of his article. Finally, I thank my girlfriend Sofia for her love and support.
%-----------------------------------------------------------------------


\begin{thebibliography}{99}
\bibitem{bobylev} A.~V.~Bobylev, A.~Palczewski, J.~Schneider: A consistency result for a discrete-velocity model of the Boltzmann equation, \emph{SIAM J.~Numer.~Anal.}, \textbf{34}~(1997), No.~5, 1865-1883.
\bibitem{kurlberg} L.~Fainsilber, P.~Kurlberg, B.~Wennberg: Lattice Points on Circles and the Discrete Velocity Model for the Boltzmann Equation.
\bibitem{hardy} G.~H.~Hardy, E.~M.~Wright: \emph{An introduction to the theory of numbers}, 5th edition, Oxford~(1979).
\bibitem{kuipers} L.~Kuipers, H.~Niederreiter: \emph{Uniform distribution of sequences}, Wiley-Interscience (1974)  
\bibitem{hecke1} E.~Hecke: \"Uber die L-funktionen und den Dirichletschen Primzahlsatz f\"ur einen beliebigen Zahlk\"orper, \emph{Nachr.~Ges.~Wiss.~G\"ottingen, Math.-Phys.~Klasse}, 1917, 299-318. 
\bibitem{hecke2} E.~Hecke: Eine neue Art von Zetafunktionen und ihre Beziehungen zur Verteilung der Primzahlen. II, \emph{Math.~Z.}, \textbf{6}~(1920), 11-51.
\bibitem{hecke3} E.~Hecke: \emph{Lectures on the Theory of Algebraic Numbers} Springer-Verlag (1981)(Translation of \emph{Vorlesung \"uber die Theorie der algebraischen Zahlen} (1923)).
\bibitem{lang} S.~Lang: \emph{Algebraic Number Theory}, Addison-Wesley (1917)
\bibitem{kubilius} J.~Kubilius: The distribution of Gaussian primes in sectors and contours, \emph{Leningrad.~U\v c.~Zap.}, \textbf{13}7~(1950), 40-52.
\bibitem{erdos} P.~Erd\H os, R.~R.~Hall: On the angular distribution of Gaussian integers with fixed norm, \emph{Disc.~Math.}, \textbf{200}~(1999), 87-94.
\bibitem{landau1} E.~Landau: \emph{Vorlesungen \"uber Zahlentheorie}, Verlag von S.~Hirzel in Leipzig (1927).
\bibitem{landau2} E.~Landau: \emph{Einf\"uhrung in die elementare und analytische Theorie der algebraischen Zahlen und der Ideale, Zweite Auflage}, Leipzig (1918), reprinted by Chelsea Publishing Company, New York (1949).
\bibitem{landau3} E.~Landau: Neuer Beweis des Primzahlsatzes und Beweis des Primidealsatzes, \emph{Math.~Ann.} \textbf{56}~(1903), 645-670.
\bibitem{landau4} E.~Landau: \emph{Handbuch der Lehre von der Verteilung der Primzahlen, Band II}, Chelsea (1953).
\bibitem{davenport} H.~Davenport: \emph{Multiplicative Number Theory. Second Edition. Revised by Hugh L. Montgomery}, Springer-Verlag (1980).
\bibitem{tenenbaum} G.~Tenenbaum: \emph{Introduction to analytic and probabilistic number theory}, Cambridge University Press (1995).
\bibitem{ireland} K.~Ireland, M.~Rosen: \emph{A Classical Introduction to Modern Number Theory, Second Edition}, Springer-Verlag (1990).
\bibitem{korner} T.~W.~K\"orner: \emph{Fourier Analysis}, Cambridge University Press (1988).
\bibitem{katai} I.~K\'atai: The distribution of divisors mod 1, \emph{Acta Math.~Hung.}, \textbf{27}~(1976), 149-152.
\bibitem{kornyei} I.~K\'atai, I.~K\"ornyei: On the distribution of lattice points on circles, \emph{Ann.~Univ.~Budapest}, \textbf{19}~(1976), 87-91.
\bibitem{vinogradov} I.~M.~Vinogradov: \emph{The method of trigonometrical sums in the theory of numbers}, Interscience Publishers (1950).
\bibitem{cilleruelo} J.~Cilleruelo: The Distribution of Lattice Points on Circles, \emph{J.~Number Theory}, \textbf{43}~(1993), 198-202.
\bibitem{turan} P.~Erd\H os, P.~Tur\'an: On a problem in the theory of uniform distribution, I\&II, \emph{Ind.~Math.}, \textbf{10}~(1948), 370-378, 406-413.   
\end{thebibliography}
\end{document}